\documentclass[msom, nonblindrev]{informs3} %
\usepackage{multibib}
\newcites{app}{References}

\usepackage{tabularx}
\newcolumntype{R}{>{\raggedleft\arraybackslash}X}
\usepackage{textcomp}
\usepackage{algorithm}
\usepackage{algpseudocode}
\usepackage{url}
\usepackage{bbm}
\usepackage{graphicx}
\usepackage{tabularx}
\usepackage{mathrsfs} 
\usepackage{amsfonts,amssymb} 
\usepackage{threeparttable}
\usepackage{amsmath}
\usepackage{color}
\usepackage{titlesec}
\usepackage{indentfirst}  
\usepackage{enumerate}
\usepackage{bm} 
\usepackage{booktabs} 
\usepackage{makecell}
\usepackage{array}
\usepackage{multirow}
\usepackage{longtable}
\usepackage{subfigure}
\usepackage{balance}
\usepackage{color}
\usepackage[colorlinks, linkcolor=blue,
            citecolor=blue ]{hyperref}
\usepackage{changes}
\usepackage{rotating}
\usepackage{appendix}
\usepackage{adjustbox}
\usepackage{lscape}
\newtheorem{Theorem}{Theorem}

\newtheorem{Lemma}{Lemma}
\newtheorem{Corollary}{Corollary}
\newtheorem{Proposition}{Proposition}
\newtheorem{Remark}{Remark}

\usepackage{natbib}
 \bibpunct[, ]{(}{)}{,}{a}{}{,}%
 \def\bibfont{\small}%

\usepackage{lineno}

\OneAndAHalfSpacedXI 
\TheoremsNumberedThrough     
\EquationsNumberedThrough    

\allowdisplaybreaks[4] 
\usepackage{setspace}

\begin{document}


\TITLE{Recommend-to-Match with Random Supply Rejections: Formulation, Approximation, and Analysis}


\ARTICLEAUTHORS{%
\AUTHOR{Haoyue Liu}
\AFF{Division of Logistics and Transportation, Shenzhen International Graduate School, Tsinghua University, Shenzhen 518055, China; \EMAIL{liu-hy22@mails.tsinghua.edu.cn} \URL{}}
\AUTHOR{Sheng Liu\footnote{This is an updated version of a paper with the same title first posted on Oct 21, 2025, on arXiv. We are grateful for the helpful comments received from Rad Niazadeh (Chicago Booth), who pointed us to a related work and motivated us to build new connections. We also thank Long He (George Washington) for the valuable input on the initial version of the paper.}}
\AFF{Rotman School of Management, University of Toronto, Toronto, Ontario M5S 3E6, Canada;\\
\EMAIL{sheng.liu@rotman.utoronto.ca}  \URL{}}
\AUTHOR{Mingyao Qi}
\AFF{Division of Logistics and Transportation, Shenzhen International Graduate School, Tsinghua University, Shenzhen 518055, China; \EMAIL{qimy@sz.tsinghua.edu.cn} \URL{}}
 } 

\ABSTRACT{\textbf{Problem definition:} Matching demand with supply in crowdsourcing logistics platforms must contend with uncertain worker participation. Motivated by this challenge, we study a two-stage “recommend-to-match” problem under stochastic supplier rejections, where each demand is initially recommended to multiple potential suppliers prior to final matching decisions. \textbf{Methodology/results:} We formulate a stochastic optimization model that explicitly captures uncertain supplier acceptance behavior. For the special case with homogeneous and independent acceptance responses, an exact mixed-integer linear program and LP formulations are achievable, but the general problem does not admit an efficient formulation. Particularly, our analysis reveals that deterministic linear approximation methods can perform arbitrarily poorly in such settings. To overcome this limitation, we propose a new approximation approach based on a convex relaxation of the original problem that admits a mixed-integer exponential cone program (MIECP) formulation. We analyze the structural properties of this approximation and establish its parametric performance guarantees. We also characterize conditions under which it can dominate a deterministic approximation. \textbf{Managerial implications:}  Extensive experiments on synthetic data and real-world freight data validate the effectiveness of our approach. Our MIECP-based solution achieves near-optimal matching performance while reducing computation time by over 90\% compared to benchmark methods, which makes it particularly promising for large-scale matching problems.
}


\maketitle


\section{Introduction}

The rise of the freelance economy has greatly revolutionized the labor market, fostering the emergence of the gig model, which promises to offer more efficient and flexible services, notably in the transport and logistics market. By 2023, the global gig economy was worth \$3.7 trillion \citep{sia}. Recent survey data suggest that more than 90\% of US workers would consider freelancing or independent contracting work \citep{GigEconomy}. Gig work encompasses a wide range of occupations, ranging from in-home care and long-haul trucking to more immediate services such as ridesharing and food delivery. Despite variations in the services offered across platforms, their core operational strategy remains the matching of customers (demand) with workers (supply).

Performing high-quality matches between demand and supply in crowdsourcing platforms poses operational challenges on multiple fronts. First, most platforms cannot dictate the choice of freelance workers, so workers may reject a match based on their individual preferences, resulting in matching inefficiency for the platform. Second, in practice, the scale of supply and demand can be large, making it challenging to compute matches within a short time frame while accounting for uncertainties. In the absence of supplier rejections, the matching problem often reduces to a linear program (LP) that scales well for large instances, which forms the foundation for developing more sophisticated dynamic matching policies. However, factoring in uncertain acceptance/rejection behaviors can break down the linearity, and the corresponding sampling-based approaches become computationally intensive as the problem size grows. 

  \begin{figure}[ht]
    \centering
    \subfigure[Order information]{
        \includegraphics[width=0.33\textwidth]{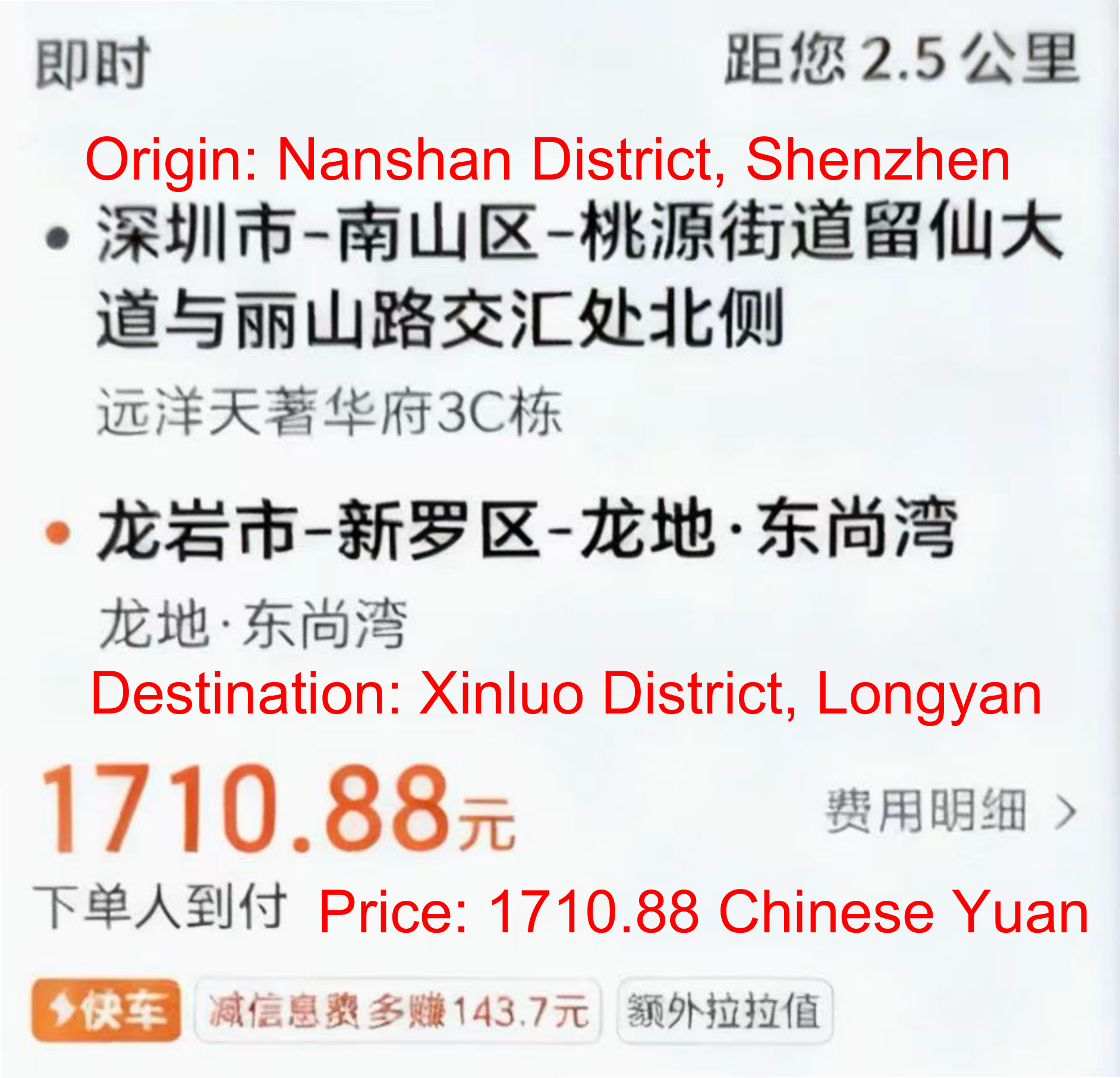}
    }
    \subfigure[Requirement and acceptance button]{
        \includegraphics[width=0.33\textwidth]{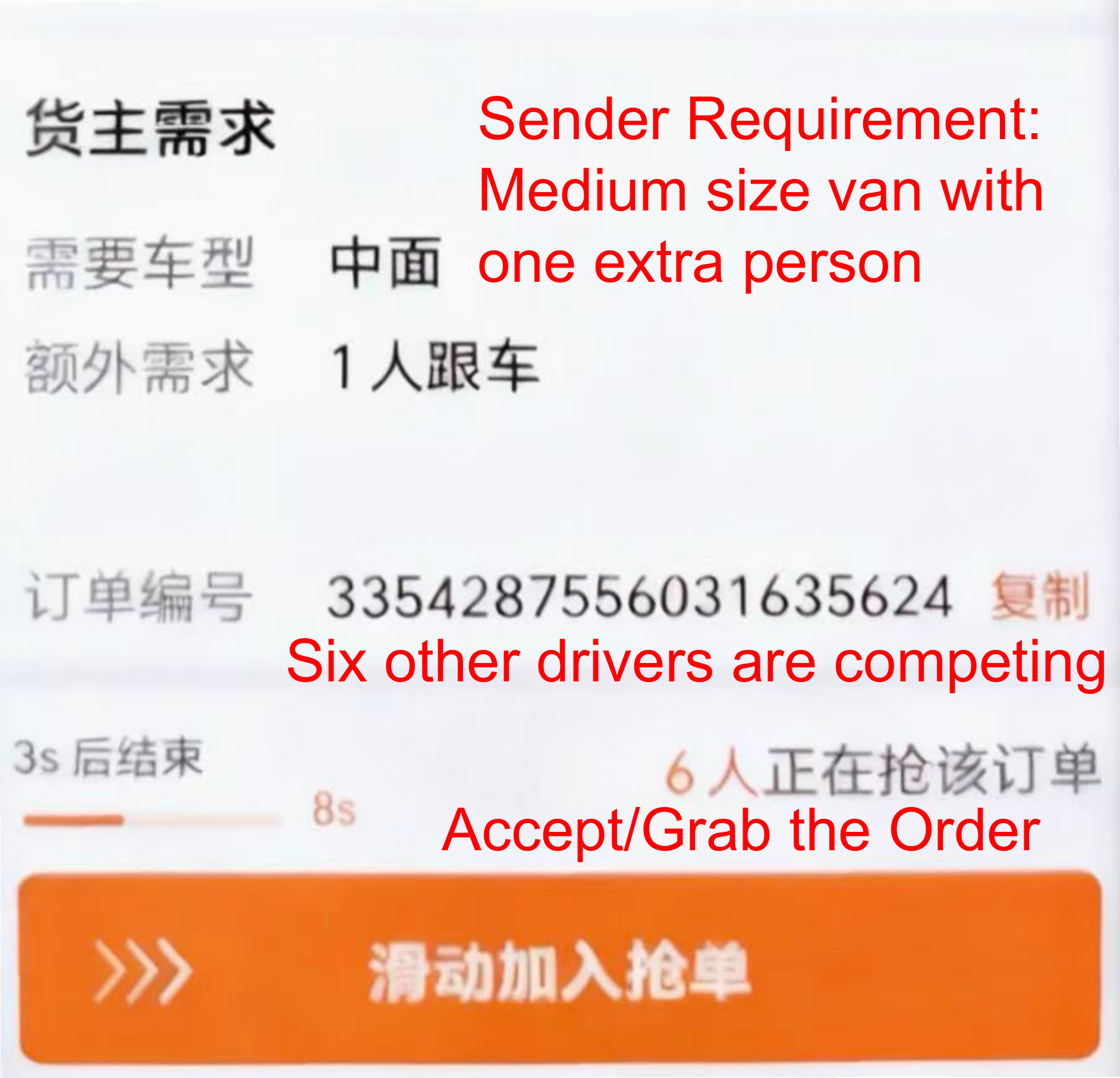}
    }
    \caption{Screenshot of the driver interface of a freight matching platform: the basic order information (origin, destination, and price) is displayed on the left, whereas the extra requirement and the acceptance button are shown on the right.}
    \label{fig:appfig}
\end{figure}

Motivated by the operations of a large freight matching platform (the screenshot of a driver interface is shown in \hyperref[fig:appfig]{Figure \ref*{fig:appfig}}) that needs to match freight and gig drivers on a daily basis (in large batches), we study a two-stage matching process called ``recommend-to-match". In this mechanism, each demand request is first recommended to at most $\theta\in \mathcal{N}^+$ suppliers. After receiving the recommended demand (match), the supplier either accepts (responds ``yes") or rejects (responds ``no") the recommendation. Then, the platform observes the supplier responses and performs a one-to-one match between demands and suppliers who respond ``yes" (commitment stage) to maximize total utility. This strategy aims to boost the likelihood of a successful final match by building a flexible buffer of suppliers, particularly when facing a large supply pool. Interestingly, ride-hailing companies have begun exploring a similar matching strategy, called non-exclusive dispatch, to broadcast ride requests to multiple drivers simultaneously before a final match is confirmed \citep{qin2025two}. In fact, Lyft has built a large-scale discrete-event simulator to test the non-exclusive dispatch policy and has found positive improvements in match time and quality \citep{ekbatani2026lyft1,ekbatani2026lyft2}. This idea is also applicable to volunteer platforms where volunteer responses are highly uncertain and depend on job attributes \citep{shi2020improving}.


The general stochastic recommend-to-match problem is challenging to solve at scale due to the nonlinearity induced by the uncertain response. The platform must decide which suppliers to recommend for each demand before observing their responses. 
Our goal in this paper is to formulate and develop efficient solution approaches to the recommend-to-match problem. \textcolor{black}{We show that LP-based matching policies using deterministic approximations can perform arbitrarily poorly and propose alternative approximation schemes with performance guarantees. The key contributions are summarized as follows.}
\begin{enumerate}
    \item Modeling. We model the recommend-to-match problem as a general stochastic optimization problem on a bipartite graph that can extend to the study of process flexibility and network design problems with arc disruptions. When the acceptance probabilities are homogeneous and independent, the problem admits a mixed-integer linear program (MILP) formulation and is further equivalent to an LP. However, the general recommend-to-match problem with heterogeneous acceptance probabilities and/or correlated responses does not yield an efficient formulation. To the best of our knowledge, we are the first to formally analyze this problem.
    \item Poor performance of direct assignment and LP-based approximation. We prove that a direct assignment policy (as common in other matching problems) with a deterministic linear approximation can perform arbitrarily poorly. Similarly, the performance of an LP-based approximation that exchanges the expectation and maximization can also be arbitrarily bad. These worst-case results caution against adopting LP-based approximation policies despite their computational tractability. 
    \item Exponential cone approximation and performance guarantee. We develop a new matching policy based on exponential cone program (ECP) that exploits the problem's structural properties. We establish parametric approximation performance guarantees under specialized and general conditions, including heterogeneous and correlated supplier acceptance. In particular, we prove that the proposed matching policy does not suffer from the poor performance of the LP-based policies, and describe conditions under which it strictly dominates deterministic LP-based approximation. 
    \item Empirical validation. We conduct extensive numerical experiments using synthetic and real-world data from a freight-matching platform to validate the superior performance of the ECP-based matching policy relative to existing benchmarks in stochastic optimization. The proposed policy achieves near-optimal expected utility while being one to two orders of magnitude faster than the conventional sampling-based approach. Compared to an alternative submodular maximization scheme, our policy delivers better solutions with significantly shorter computational time. The advantage of the ECP-based policy is particularly pronounced for large-scale applications.
\end{enumerate}

\section{Literature review}
We review three streams of related literature: the general two-sided matching problem (TSMP), matching in transportation and logistics, and the approximation of stochastic matching-related problems. 

\subsection{Two-sided Matching Problems}

According to matching cardinality constraints, two-sided matching problems are typically categorized into three types: one-to-one, many-to-one (also called one-to-many), and many-to-many matching \citep{gale1962college, roth1984stability, echenique2006theory}. Our study focuses on a specific one-to-one matching, where each demand-side participant can be recommended to multiple supply-side participants, each supply-side participant can be recommended to at most one demand-side participant, corresponding to a many-to-one recommendation, and the final matching is one-to-one. For a comprehensive overview of TSMPs, we refer the reader to \citet{roth1992two} and \citet{abraham2003algorithmics}. Below, we will mainly focus on reviewing TSMPs related to our research that deal with heterogeneity and uncertainty on both sides.

The existing research primarily focuses on the sequential matching decisions without rejections \citep{ashlagi2019matching,baccara2020optimal,hu2022dynamic,kerimov2024dynamic,kerimov2025optimality} or matching with preference uncertainties 
\citep{aziz2022stable,stokkink2025optimal}. Related to our setting, recent studies on two-sided assortment optimization consider platforms that influence matching outcomes through recommendation sets under decentralized choice on both sides \citep{ashlagi2022assortment,aouad2023online}. In particular, \citet{ashlagi2022assortment} study a two-sided sequential assortment model with Multinomial Logit (MNL)-type customer and supplier choices, where the platform optimizes the menu profile to maximize expected matches, whereas \citet{aouad2023online} analyze online assortment policies under dynamic arrivals and supplier-side MNL and Nested Logit models. However, these papers focus on menu design and online assortment policies rather than explicit recommend-to-match optimization with stochastic supplier acceptance and a subsequent platform-controlled commitment stage.  
A few recent papers have looked more directly into the uncertainty of agent abandonment/rejections and incorporated it into matching optimization.  
For example, \citet{hou2025reinforced} address matching instability in crowdsourced delivery caused by stochastic driver acceptance behavior and propose a reinforced stable matching mechanism, which uses the Gale-Shapley algorithm to find an initial stable match, and then employs a stochastic linear program model to determine an optimal compensation to incentivize drivers. 
\citet{wang2025order}’s study is most closely related to ours as they analyze cancellation behavior following customer–driver recommendations in ride-sourcing markets and maximize platform revenue under different pricing schemes. Nevertheless, their model restricts each demand to be recommended to at most one supply at the recommendation stage, which could result in a lower matching success rate than what is achievable under our schemes.


\subsection{Matching in Transportation and Logistics}
Matching constitutes a core decision in transportation and logistics services, such as vehicle–cargo matching in long-haul transportation and driver–order matching in on-demand delivery \citep{liu2021time, zhao2024market}. 
One line of research focuses on preference estimation of couriers and then solving a bipartite matching model based on revenue maximization. \citet{deng2021prediction} introduce a framework to predict the probability of vehicle–cargo matching using a dynamic Bayesian network, where the matching success rate and costs vary over time. 
\citet{wang2022recommending} and \citet{wang2023online} develop a tree-based model and a deep reinforcement learning model, respectively, to predict the driver's acceptance behavior and preference in food delivery systems. These studies focus on predictive modeling instead of detailed matching optimization.

Recent studies on crowdsourcing logistics have explored advanced matching methods to cope with operational challenges, which can be categorized into four types according to \cite{mohri2023crowdshipping}: menu offering \citep{mofidi2019beneficial}, in-store customers \citep{dayarian2020crowdshipping}, alternate delivery points \citep{kizil2023public}, and employing transshipment points \citep{kafle2017design}. Among these, the menu offering approach is most relevant to our study, as it recommends multiple drivers for each order and allows drivers to self-select orders from personalized menus designed by the platform. \citet{mofidi2019beneficial} first explore the single-stage deterministic menu offering problem, which they formulate as a mixed-integer bilevel program, and introduce a specialized reformulation technique. \citet{horner2021optimizing} then extend \citet{mofidi2019beneficial}'s work to a stochastic setting that accounts for uncertainty in driver-side utilities and relaxes the constraints on the number of orders that a driver can select and serve. They reformulate the original bilevel model into a single-level formulation and address uncertainty using the Sample Average Approximation (SAA) method. \citet{ausseil2022supplier} investigate this approach in a multi-stage dynamic setting, where the subproblem in each stage is similar to \citet{mofidi2019beneficial}'s formulation. \citet{karabulut2022value} adopt the model used by \citet{mofidi2019beneficial} and introduce preference learning module in optimization model. Specifically, they iteratively optimize, observe supplier choices, learn preferences, and refine the utility estimation function in the optimization model in the dynamic multi-period system. More recently, \citet{horner2025increasing} incorporates decisions of driver-side pricing into the menu offering approach, assuming a linear relationship between driver acceptance probability and the offered price. They formulate a SAA linear integer model and show that personalized pricing can increase drivers’ willingness to accept requests, benefiting both supply and demand sides. 


Our recommend-to-match problem is distinct from the menu-offering studies, both practically and methodologically. From the application perspective, these studies are set in ridesharing and crowdsourced delivery platforms, where demand and supply sizes are relatively comparable \citep{karabulut2022value}. In such contexts, menu offerings help increase demand acceptance rates, but designing customized menus for every driver is cumbersome, especially when the supply pool size is large. Motivated by the freight matching platform, where the supply pool can be larger than the demand pool, the recommend-to-match mechanism enables the platform to exercise greater control over final matches. Methodologically, stochastic models in the menu offering literature often rely on SAA to address uncertainty. However, SAA-based methods may not scale well for large instances, and the existing studies rarely tackle instances with more than 50 agents per side \citep{horner2021optimizing, ausseil2022supplier}. In contrast, we design a scalable ECP-based approximation that does not rely on sampling while still enjoying performance guarantees. \textcolor{black}{Our recommend-to-match mechanism also differs from the one-to-many match policies developed for community first-responder and emergency dispatch systems, which adopt a first-accept rule, assigning the job to the first timely responder to ensure rapid response \citep{henderson2022should,dellaert2024community}.}

The closest to our studied problem in transportation applications is \cite{ekbatani2026lyft1}, where the authors independently study a recommend-to-match problem for ride-hailing. \textcolor{black}{Specifically, they frame the problem as a notification set selection problem and study both a first-acceptance mechanism and a best-acceptance mechanism in the second stage, with the latter being more closely related to our setting.}
Technically, for the best-acceptance setting, they exploit the structure of the objective function and develop approximation algorithms using a submodular optimization oracle, which holds when the drivers' response behaviors are independent. In contrast, our solution framework is based on a convex relaxation of the problem with a tuning parameter that can yield a tighter approximation and apply to more general scenarios with correlated responses. 

\subsection{Approximation of Stochastic Matching-related Problems}


Another branch of related research focuses on developing efficient algorithms with performance guarantees for stochastic matching, where some nodes or edges may not be present or might fail with a certain probability. The main approximation approaches can be categorized into non-adaptive and adaptive methods. \citet{chen2009approximating} originally study the unweighted stochastic matching problem and provide a 4-approximation with the non-adaptive greedy algorithm, which is later proved to be a 2-approximation by \citet{adamczyk2011improved}. However, the greedy approach performs poorly for the weighted setting, motivating subsequent work on adaptive methods that sequentially probe edges according to priority, achieving better approximation guarantees. \citet{bansal2012lp} consider the weighted stochastic matching problem and provide a 3-approximation for bipartite graphs and a 4-approximation for general graphs using adaptive LP rounding techniques. For stochastic matching with edge existence probability and nodes' timeouts (the upper-bound number of edges incident to a node that can be probed), \citet{adamczyk2015improved} and \citet{baveja2018improved} introduce non-uniform random ordering and LP-guided probing probabilities to achieve constant approximation ratios. Recently, \citet{you2024approximate} study a general class of dynamic stochastic matching problems using approximate dynamic programming and establish an upper bound on the value of the optimal policy. 

Although our problem shares certain features of the classical stochastic matching problem, such as edges appearing with probabilities and vertices having probing limits, the objective function differs in that the value of an edge depends on the recommendation decision rather than being independent, which is a core assumption in previous research that guarantees the approximation bound \citep{adamczyk2015improved, baveja2018improved}. Because each demand can be recommended to multiple suppliers with heterogeneous acceptance probabilities in our setting, the objective function becomes implicitly nonlinear, and LP-guided approaches can no longer yield effective approximations.

As a final note, although the exponential cone program has been proposed for other problems such as product design and scheduling \citep{akccakucs2021exact, chen2024exponential}, we are the first to develop cone program-based matching policies. Our performance analysis leverages unique features from supply rejections and matching conditions that could be useful for developing effective dynamic matching policies.


\section{Model and Analysis}\label{sec:problem}
In this section, we first introduce the recommend-to-match problem, followed by exact formulations. Then we discuss and develop approximation policies and analyze their performance. We use supply and supplier interchangeably whenever no confusion may arise.
 
\subsection{Problem Description and Formulation}

We consider a one-shot matching problem where supply has a probability of rejecting the demand, resulting in a match failure. To enhance the likelihood of successful matches, an effective approach is to introduce a recommendation stage, in which each demand is recommended to multiple suppliers (while each supplier is recommended at most one demand each time). The system then observes supplier choices and subsequently conducts the final matching between demand and supply.\footnote{The recommend-to-match approach introduces demand flexibility but limits the supply choice to hedge against uncertain supply behavior that is complicating the platform operations, especially when the supply pool is large and heterogeneous.} In this problem, the finite sets of supplies and demands are denoted as $\mathcal{A}^{S}$ and $\mathcal{A}^{D}$, respectively. Only matching between demand $i \in \mathcal{A}^{D}$ and supply $j \in \mathcal{A}^{S}$ pairs is feasible, each corresponding to a positive utility $u_{ij}$ and acceptance probability $p_{ij}$. For each demand, it can be recommended to at most $\theta$ supplies and finally assigned to the one with the highest utility among the accepted supplies. The decision variable $x_{ij}$ indicates whether to recommend demand $i$ to supply $j$. To capture the uncertainties of supply acceptance, we define random variables $\tilde{\xi}_{ij}$ that follow a Bernoulli distribution with parameter $p_{ij}$: $\tilde{\xi}_{ij}=1$ if supply $j$ accepts demand $i$, and 0 otherwise. The acceptance uncertainties are independent across supplies, but there could be correlation among demands (we will discuss relaxing this assumption later in Section \ref{sec:analysis}). The joint distribution of $\tilde{\xi}_{ij}$ for all $i \in \mathcal{A}^{D}$ and $j\in \mathcal{A}^{S}$ is denoted as $\mathbb{P}$. We denote the supply-to-demand ratio $\frac{|\mathcal{A}^{S}|}{|\mathcal{A}^{D}|}$ as $\gamma$, whose round-down and round-up are $\lfloor \gamma \rfloor$ and $\lceil \gamma \rceil$, respectively. The recommend-to-match problem is then formulated as the following stochastic program:
\begin{subequations}
    \begin{align}
        \textbf{[SP]}:  \max \ & \sum_{i \in \mathcal{A}^{D}} \mathbb{E}_{\mathbb{P}} \left [\max_{j\in \mathcal{A}^{S}} u_{ij} \tilde{\xi}_{ij} x_{ij}\right] \label{eq:a0}\\
        {\rm s.t.} \ & \sum_{j\in \mathcal{A}^{S}}x_{ij}  \le \theta , \forall i \in \mathcal{A}^{D}, \label{eq:a1}\\
        &\sum_{i \in \mathcal{A}^{D}}x_{ij} \le 1, \forall j\in \mathcal{A}^{S},\label{eq:a2}\\
        & x_{ij} \in \{0,1\}, \forall i \in \mathcal{A}^{D}, j\in \mathcal{A}^{S}\label{eq:a3},
    \end{align}
\end{subequations}
where the objective function \hyperref[eq:a0]{(\ref*{eq:a0})} maximizes the expected total utility, constraints \hyperref[eq:a1]{(\ref*{eq:a1})} restrict each demand to be recommended to a maximum of $\theta$ supplies, constraints \hyperref[eq:a2]{(\ref*{eq:a2})} guarantee that each supply is recommended to at most one demand, and constraints \hyperref[eq:a3]{(\ref*{eq:a3})} define the domain of decision variables in this model. The main difficulty in solving SP lies in evaluating the objective function, which involves an inner maximization perplexed by the Bernoulli random variables.\footnote{There are potentially other operations constraints facing the platform, but we keep the model parsimonious here to highlight the key tradeoffs. } 

\begin{Remark}
    The structure of the above formulation is not unique to the recommend-to-match problem, but also appears in other operations applications. For example, in process flexibility design, $x_{ij}$ corresponds to the link decision between a plant/manufacturer and a product, $u_{ij}$ represents the unit profit from production, and $\tilde{\xi}_{ij}$ indicates the disruption uncertainty of the link. In network design problems, $x_{ij}$ may indicate the assignment of customers to (primary and backup) facilities, in addition to the facility-opening decisions; when random transportation blockages occur, reassignment can be made to save costs. This way, SP could also model a facility location or network design problem with arc-disruption risks \citep{mak2009stochastic,mehmanchi2020analysis,salman2015emergency,hassin2017multiple}. 
\end{Remark}

We first consider a special case where the acceptance probabilities $p_{ij}=p$ for all demand-supply pairs, i.e., the acceptance probabilities are homogeneous. Furthermore, assuming the acceptance behaviors are independent across supplies, then the probability of assigning a demand to the supply with $r$th highest utility ($r\in \mathcal{R}$) among the recommended supplies is $p(1-p)^{r-1}$, which we denote as $p_r$. In this case, SP can be reformulated as the following MILP:
\begin{subequations}
    \begin{align}
        \textbf{[R-SP]}:  \max \ & \sum_{i \in \mathcal{A}^{D}}\sum_{r \in \mathcal{R}}p_{r}w_{ir} \label{eq:r0}\\
        {\rm s.t.} \ & \hyperref[eq:a1]{(\ref*{eq:a1})}, \hyperref[eq:a2]{(\ref*{eq:a2})}, \nonumber  \\
        & z_{ijr} \le x_{ij}, \forall i \in \mathcal{A}^{D}, j\in \mathcal{A}^{S}, r \in \mathcal{R}, \label{eq:r2}\\
        &  t_{ir} \le r, \forall i \in \mathcal{A}^{D}, r \in \mathcal{R},\label{eq:r4}\\
        &  t_{ir} \ge r-\theta(1-v_{ir}), \forall i \in \mathcal{A}^{D}, r \in \mathcal{R},\label{eq:r5}\\
        &  t_{ir} \ge \sum_{j\in \mathcal{A}^{S}}x_{ij}-\theta v_{ir}, \forall i \in \mathcal{A}^{D}, r \in \mathcal{R},\label{eq:r6}\\
         & \sum_{j\in \mathcal{A}^{S}}z_{ijr}=t_{ir} , \forall i \in \mathcal{A}^{D}, r \in \mathcal{R},\label{eq:r7}\\
         & w_{ir} \le u_{ij}+ M^{1}_{i}(1-z_{ijr}),\forall i \in \mathcal{A}^{D}, j\in \mathcal{A}^{S}, r \in \mathcal{R},\label{eq:r8}\\
         & w_{ir} \le M^{1}_{i}v_{ir}, \forall i \in \mathcal{A}^{D}, r \in \mathcal{R},\label{eq:r9}\\
         & x_{ij},z_{ijr},v_{ir} \in \{0,1\}, t_{ir} \in \mathbb{N}, w_{ir} \ge 0, \forall i \in \mathcal{A}^{D}, j\in \mathcal{A}^{S},  r \in \mathcal{R}.\label{eq:r10}
    \end{align}
\end{subequations}

The objective function \hyperref[eq:r0]{(\ref*{eq:r0})} maximizes the weighted sum of recommended supplies' utility from the first to the $\theta$-th largest for each demand. Besides the original constraints in SP model, the R-SP model introduces other constraints to obtain the utility value of the $r$ largest recommended supplies for demand $i$, which is denoted as $w_{ir}$. 
Constraints \hyperref[eq:r2]{(\ref*{eq:r2})} indicate that supply $j$ can only be among the top $r$ utility supplies for demand $i$ if the demand-supply pair $(i,j)$ is recommended. Constraints \hyperref[eq:r4]{(\ref*{eq:r4})-(\ref*{eq:r6})} linearize the relationship that $t_{ir}$ is the minimum value of $r$ and $\sum_{j\in \mathcal{A}^{S}}x_{ij}$, by introducing a binary auxiliary variable $v_{ir}$. Variable $v_{ir}$ equaling 1 indicates the value of $\sum_{j\in \mathcal{A}^{S}}x_{ij}$ is not less than $r$. Constraints \hyperref[eq:r7]{(\ref*{eq:r7})} guarantee that the number of recommended supplies of demand $i$ in the top $r$ utility supplies is equal to the minimum value of $\sum_{j\in \mathcal{A}^{S}}x_{ij}$ and $r$. Constraints \hyperref[eq:r8]{(\ref*{eq:r8})} ensure if demand-supply pair $(i,j)$ is one of the top $r$ utilities for demand $i$, then the utility $w_{ir}$ must not be greater than utility $u_{ij}$. Constraints \hyperref[eq:r9]{(\ref*{eq:r9})} restrict value $w_{ir}$ equaling 0 when the value of $r$ is greater than the number of recommended supplies for demand $i$. Note that parameter $M^{1}_{i}$ in \hyperref[eq:r8]{(\ref*{eq:r8})-(\ref*{eq:r9})} is a big number, which can be set as the maximum utility value for demand $i$. The full notation of this formulation is summarized in \hyperref[tab:notation]{Table \ref*{tab:notation}}. \hyperref[pro1]{Lemma \ref*{pro1}} shows that R-SP model is a valid MILP formulation of SP when the acceptance random variables are independent and identically distributed (i.i.d.) and \hyperref[corollary0]{Corollary \ref*{corollary0}} shows that the R-SP can be reduced to an LP. 
\begin{table}[htbp]
\caption{Notation for R-SP}
\label{tab:notation}
\centering

\renewcommand{\arraystretch}{1.15}
\setlength{\tabcolsep}{8pt}

\begin{tabularx}{\textwidth}{l X}
\hline
\multicolumn{2}{l}{\textbf{Sets}}\\ \hline
$\mathcal{A}^{D}$ & Set of demands, indexed by $i$.\\
$\mathcal{A}^{S}$ & Set of supplies, indexed by $j$.\\
$\mathcal{R}$ & Set of positive integers from 1 to the maximum recommendation number, indexed by $r$.\\ \hline
\multicolumn{2}{l}{\textbf{Parameters}}\\ \hline
$\theta$ & Maximum number of supplies that can be recommended to each demand.\\
$u_{ij}$ & Platform's utility for supply $j$ fulfilling demand $i$.\\
$p$ & Homogeneous acceptance probability of supplier.\\
$p_{r}$ & The probability of assigning to the supply with $r$th highest utility among the recommended suppliers.\\
$\gamma$ & Ratio of $|\mathcal{A}^{S}|$ and $|\mathcal{A}^{D}|$; its round-down and round-up are denoted by $\lfloor\gamma\rfloor$ and $\lceil\gamma\rceil$, respectively.\\ \hline
\multicolumn{2}{l}{\textbf{Decision Variables}}\\ \hline
$x_{ij}$ & Binary variable, taking 1 if the platform recommends demand $i$ to supply $j$, 0 otherwise.\\
$z_{ijr}$ & Binary variable, taking 1 if supply $j$ is among the top $r$ utility supplies in demand $i$'s recommendations, 0 otherwise.\\
$t_{ir}$ & Integer variable, indicating the minimum of $\sum_{j\in \mathcal{A}^{S}}x_{ij}$ and $r$.\\
$v_{ir}$ & Binary variable, taking 1 if $\sum_{j\in \mathcal{A}^{S}}x_{ij} \ge r$, 0 otherwise.\\
$w_{ir}$ & Continuous variable, indicating the utility of the $r$th largest recommended supply in demand $i$, set to 0 if it does not exist.\\
\hline
\end{tabularx}
\end{table}
\begin{Lemma}\label{pro1}
    When acceptance random variables $\tilde{\xi}_{ij}$ are i.i.d. across $i$ and $j$, the optimal objective value and the set of optimal solutions of model SP coincide with those of R-SP.
\end{Lemma}

\begin{Corollary} \label{corollary0} The R-SP formulation can be reduced to an LP as described in \citet{ekbatani2026lyft1}.   
\end{Corollary}

We note that the proof of \hyperref[corollary0]{Corollary \ref*{corollary0}} is based on reformulating the R-SP model, of which the core lies in the monotonicity property of the rank index $r$ for utility-related variables ($w_{ir}$ in MILP). To link the MILP and the LP formulation, we augment the MILP with a new formulation and show that the constraints associated with the binary variables do not tighten the feasible region and are therefore redundant, leading to an equivalent LP formulation presented for a similar problem in \cite{ekbatani2026lyft1}.

For general cases where the probabilities $p_{ij}$ are heterogeneous or there is correlation among supply acceptances, SP cannot be reformulated to a tractable MILP, as we have to track the realization of acceptance of each supply individually. One alternative is to employ an SAA scheme. Let $\mathcal{S}$ denote the set of samples drawn from the acceptance distribution and $\xi_{ij}^{s}$ be the acceptance response of supply $j$ for demand $i$ in sample $s \in \mathcal{S}$. We introduce a binary variable $\omega_{ij}^{s}$, indicating whether demand $i$ is finally assigned to supply $j$ in sample $s$. Then, the SAA formulation for the recommend-to-match problem is given by
\begin{subequations}
    \begin{align}
        \textbf{[SP-SAA]}:  \max \ & \frac{1}{|\mathcal{S}|}\sum_{s \in \mathcal{S}}\sum_{i \in \mathcal{A}^{D}}\sum_{j\in \mathcal{A}^{S}}u_{ij}\omega_{ij}^{s}  \label{eq:b0}\\
        {\rm s.t.} \ & \hyperref[eq:a1]{(\ref*{eq:a1})}, \hyperref[eq:a2]{(\ref*{eq:a2})}, \nonumber \\
        &\sum_{j\in \mathcal{A}^{S}} \omega_{ij}^{s} \le 1, \forall i \in \mathcal{A}^{D}, s \in \mathcal{S} \label{eq:b2}\\
        &  \omega_{ij}^{s} \le \xi_{ij}^{s}x_{ij},  \forall i \in \mathcal{A}^{D}, j\in \mathcal{A}^{S}, s \in \mathcal{S} \label{eq:b3}\\
        & x_{ij} \in \{0,1\}, \omega_{ij}^{s} \in \{0,1\}, \forall i \in \mathcal{A}^{D}, j\in \mathcal{A}^{S},  s \in \mathcal{S}\label{eq:b4} .
    \end{align}
\end{subequations}

SP-SAA maximizes the sample average total utility while satisfying the constraints of recommendation and supply selection. Constraints \hyperref[eq:b2]{(\ref*{eq:b2})} ensure each demand is finally assigned to one supply in every sample. Constraints \hyperref[eq:b3]{(\ref*{eq:b3})} indicate that a demand-supply pair in a sample can only be matched when it is recommended and the supply is accepted. Although SP-SAA can be solved by off-the-shelf solvers, the solution time can grow substantially as the sample size increases for large-scale instances (the number of binary variables is $|\mathcal{A}^D|\cdot|\mathcal{A}^S|\cdot |S| + |\mathcal{A}^D|\cdot|\mathcal{A}^S|)$, which can prohibit its effective deployment in practice.  Therefore, our goal is instead to seek high-quality approximation algorithms that are computationally efficient. 

\subsection{Deterministic Approximation Benchmarks}
In this subsection, we analyze two natural deterministic approximation policies that admit LP-based reformulations. 
\subsubsection{Direct Assignment Approximation.}
In the spirit of fluid approximation, one may replace the inner maximization by the summation of utilities weighted by acceptance probabilities. Doing so mimics the direct assignment policy (DAP) in bipartite matching problems, where the recommendation stage is effectively diminished. Specifically, DAP corresponds to the following integer linear program, \textcolor{black}{which can be reduced to an LP because its constraint coefficient matrix is totally unimodular.}
\begin{subequations}
    \begin{align}
       \textbf{[DAP]}:  \max \ & \sum_{i \in \mathcal{A}^{D}} \sum_{j \in \mathcal{A}^{S}} p_{ij}u_{ij}x_{ij} \\
       {\rm s.t.} \ & \hyperref[eq:a1]{(\ref*{eq:a1})-(\ref*{eq:a3})} \nonumber.
    \end{align}
\end{subequations}

Let $\mathbf{x}^{\star}$ and $\mathbf{x}^{D}$ denote the optimal solution to SP and DAP, respectively. We define the approximation gap of DAP as $\frac{Obj^{E} (\mathbf{x}^{\star})-Obj^{E} (\mathbf{x}^{D})}{Obj^{E} (\mathbf{x}^{\star})}$, where $Obj^{E}(\mathbf{x}^{\star})$ and $Obj^{E}(\mathbf{x}^{D})$ denote the objective value of SP at solutions $\mathbf{x}^{\star}$ and $\mathbf{x}^{D}$, respectively. Despite the computational benefits, the following result illustrates that DAP's approximation performance is not bounded. 
\begin{Proposition} \label{prop:unbound}
    Assume that the utility of each demand-supply pair lies within the interval $[ a, b ]$ with homogeneous acceptance probabilities $p$, there exists an instance such that the approximation gap of DAP can be arbitrarily close to one, i.e., 
    \begin{align*}
    \frac{Obj^{E} (\mathbf{x}^{\star})-Obj^{E} (\mathbf{x}^{D})}{Obj^{E} (\mathbf{x}^{\star})} \rightarrow 1.
\end{align*}
\end{Proposition}

\vspace{3pt}

\hyperref[prop:unbound]{Proposition \ref*{prop:unbound}} implies that the DAP approximation can be arbitrarily bad, especially for cases with relatively large values of $\theta$ (or small values of $\gamma$). The main reason is that DAP greedily recommends a subset of demands with the highest utilities, which can result in under-utilization of supplies. This under-utilization could lead to greater losses when the supply pool is more limited or when each demand can access more supplies (leading DAP to recommend fewer demands). 

\subsubsection{Deterministic Inner Maximization Approximation.}
To address the limitation that the DAP model may recommend suppliers to only a small subset of demands, an alternative approximation strategy is to exchange the expectation and the inner maximization over suppliers, which yields the deterministic inner maximization policy (DIMP), corresponding to the following formulation:
\begin{subequations} \label{dapnew}
    \begin{align}
       \textbf{[DIMP]}:  \max \ & \sum_{i \in \mathcal{A}^{D}} \max_{j \in \mathcal{A}^{S}} p_{ij}u_{ij}x_{ij} \label{obj:dapnew} \\
      {\rm s.t.} \ & \hyperref[eq:a1]{(\ref*{eq:a1})-(\ref*{eq:a3})}  \nonumber.
    \end{align}
\end{subequations}
The objective function \hyperref[obj:dapnew]{(\ref*{obj:dapnew})} can be linearized by (i) replacing inner maximization item $\max_{j \in \mathcal{A}^{S}} p_{ij}u_{ij}x_{ij}$ with a continuous variable $v_i$, and (ii) adding a binary auxiliary variable $z_{ij}$ indicating whether the supply $j$ is the one that attains the maximum $p_{ij}\cdot u_{ij}$ among all suppliers recommended to demand $i$. \textcolor{black}{By introducing related logic constraints with the standard Big-M linearization method, DIMP can be represented as an MILP. We further observe that in the DIMP model, only the supplier $j$ with the largest utility-weighted acceptance probability contributes to the objective value. This indicates that DIMP can be equivalently reformulated as an LP by replacing $\theta$ with one in the recommendation constraints. }
However, this equivalence implies a critical weakness: the model effectively recommends only one supplier per demand, making the optimal solution independent of the recommendation number $\theta$. We now stress the weakness of DIMP with a worst-case performance analysis similar to that of DAP. Specifically, let $\mathbf{x}^{M}$ denote the optimal solution to DIMP and define the approximation gap of DIMP as $\frac{Obj^{E} (\mathbf{x}^{\star})-Obj^{E} (\mathbf{x}^{M})}{Obj^{E} (\mathbf{x}^{\star})}$. The following result shows that DIMP's approximation performance is also not bounded. 

\begin{Proposition}  \label{dapnew}Assume that $\gamma \ge \theta >1$, the utility of each demand-supply pair lies within the interval $[ a, b ]$ with homogeneous acceptance probabilities $p$, there exists an instance such that the approximation gap of DIMP can be arbitrarily close to one, i.e., 
\begin{align*}
    \frac{Obj^{E} (\mathbf{x}^{\star})-Obj^{E} (\mathbf{x}^{M})}{Obj^{E} (\mathbf{x}^{\star})} \rightarrow 1. 
\end{align*}
\end{Proposition}

Specifically, when $b\rightarrow a$, $p \rightarrow 0$, $\gamma\ge \theta$ and $\theta \rightarrow \infty$, DIMP will perform arbitrarily poorly. The main reason is exactly that, in the DIMP model, only one supply’s utility contributes to the match, pushing the total utility toward zero when each demand is recommended to only the highest-utility supplier, and the acceptance probability is very low. Meanwhile, when $\theta$ is sufficiently large, the optimal solution to the SP model can effectively mitigate this issue by connecting demand to many suppliers (e.g., diversification). Moreover, when utility values are relatively concentrated, the objective value is primarily determined by the overall acceptance probability, thereby widening the performance gap between DIMP and SP.
Although the condition stated above appears stringent, we prove that the ECP approximation proposed later dominates DIMP in more general settings, as shown in \hyperref[dominate]{Proposition \ref*{dominate}} of Section \ref{sec:analysis}. Next, we introduce a new approximation approach based on exponential cone programming that addresses the limitations of these deterministic methods.

\subsection{MIECP Approximation}\label{subsec:ecp}
We first recall a well-known property of the log-sum-exp function in \hyperref[lemma0]{Lemma \ref*{lemma0}}, which will be useful in approximating SP's inner non-linear maximum function with a performance guarantee.


\begin{Lemma}\label{lemma0}
    \citet[Section~3.1.5]{boyd2004convex} For any $\{z_i\}_{i \in [N]}$ and $\tau>0$, the following inequality holds:
\begin{align}
    \max \{z_1,z_2,...,z_n \} \le \tau \log\left(\sum_{i \in N}\exp{\left(\frac{z_i}{\tau}\right)}\right) \le\max \{z_1,z_2,...,z_n \}+\tau\log n \label{eq:property}
\end{align}
\noindent where $\tau$ is a parameter that can control the approximation accuracy.
\end{Lemma}

Based on the above lemma, we show that the objective function of SP can be bounded as follows.

\begin{Corollary} \label{corollary1} For any $\tau>0$, the objective function of SP satisfies
\[  Obj^E(\mathbf{x}) \leq  \tau\sum_{i \in \mathcal{A}^{D}}  \log \left[\sum_{j\in \mathcal{A}^{S}}\left((\exp{\left(\frac{u_{ij} }{\tau}\right)}-1)p_{ij} + 1\right)x_{ij}\right]      \]
\end{Corollary}

This bound offers a way to approximate the expectation by a nonlinear function that is more amenable to optimization. {\color{black} Importantly, this upper bound does not rely on the independence of acceptance random variables $\mathbf{\xi}$ and holds for cases with correlated driver acceptances.}

By replacing the objective function of SP with the upper bound in \hyperref[corollary1]{Corollary \ref*{corollary1}}, we obtain a convex relaxation of the original problem. A practical question that underpins this relaxation is the choice of $\tau$. To that end, we first provide a sufficient condition for the existence of optimal $\tau$, based on \hyperref[assu1]{Assumption \ref*{assu1}}, that minimizes the approximation objective function $f(\tau)=\tau \sum_{i}\log \left(\sum_{j}\Bigl(p_{ij}\exp({u_{ij}/\tau})+ (1-p_{ij})\Bigr)\right)x_{ij}$ for tightness. 
\begin{assumption}\label{assu1}
    For each demand $i\in \mathcal{A}^D$, the utilities of different suppliers are strictly different, and the acceptance probability $p_{ij}$ is strictly less than 1.
\end{assumption}
 \begin{Proposition}
\label{prop:not_infty_minimizer} Suppose each demand is recommended to at least one supplier, and there exists at least one demand that is recommended to two (or more) suppliers, i.e.,
\[
|A_i|\ge 1\ \ \forall i\in \mathcal{A}^D\quad\text{and}\quad \exists\, i_0\ \text{s.t.}\ |A_{i_0}|\ge 2,
\]
\noindent where $A_i:=\{j:\ x_{ij}=1\}$, then the minimizer of $f(\tau)$ over $\tau>0$ exists. 
\end{Proposition}

We note that \hyperref[assu1]{Assumption \ref*{assu1}} is not stringent because suppliers' preferences (which affect utilities) depend not only on driver and job attributes but also on drivers' past interaction history, which is highly variable. The sufficient conditions laid out in \hyperref[prop:not_infty_minimizer]{Proposition \ref*{prop:not_infty_minimizer}} are likely to hold when the supplier pool is sufficiently large. Despite its existence, the optimal $\tau$ does not admit a closed-form expression, and more critically, its value depends on the recommendation decisions $x_{ij}$. This implies that there is no universal optimal $\tau$ that yields the tightest bound. Additionally, the tightness of the bound may not necessarily translate to the superiority of the approximation solution, which will be investigated next. 

\subsubsection{Choice of $\tau$.}\label{sec:choice} To understand how the choice of $\tau$ affects the solution performance, we conduct a comprehensive empirical investigation over instances of different sizes, where we numerically search for the optimal $\tau$ that minimizes $f(\tau)$ on 100 randomly generated solutions ($x$) for each instance. Specifically, we consider instances denoted by D$i$-S$j$ with $i$ demands and $j$ suppliers. These ten types of instances range from small to moderate scale and include $(i,j)\in\{(5,10), (5,15), (5,20), (5,25), (10,20), (10,30), (10,40), (10,50), (20,50), (20,80)\}$. Other parameters are set as follows: maximum recommendation size $\theta \in \{2,3,4,5\}$, acceptance-probability interval in $\{(0.7,0.75), (0.7,0.8), (0.7,0.85), (0.7,0.9)\}$, and utility values drawn from the interval $[0.4,1]$ uniformly. Summing up, 160 instance types are created, and for each type, 50 random instances are generated, giving a total of 8,000 test instances. 

\begin{figure}[htbp]
    \centering
    \subfigure[Optimal $\tau$ that minimizes the value of $f(\tau)$.]{
        \includegraphics[width=0.65\textwidth,trim=20 0 20 35, clip]{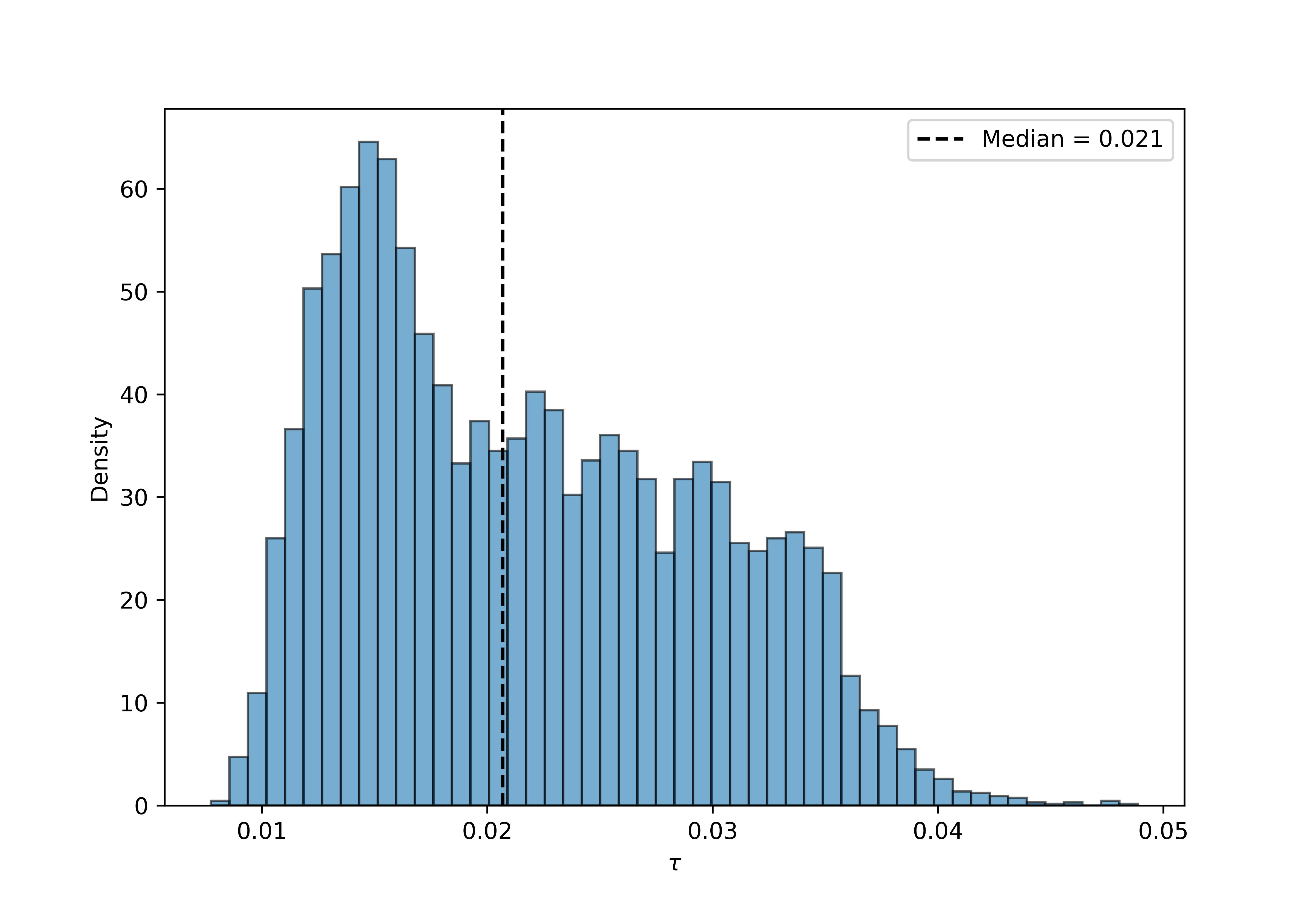}
        \label{fig:tauC_1}
    }
    \hfill
    \subfigure[Optimality gap of different $\tau$ 
    ]{
        \includegraphics[width=0.65\textwidth,trim=20 0 20 35, clip]{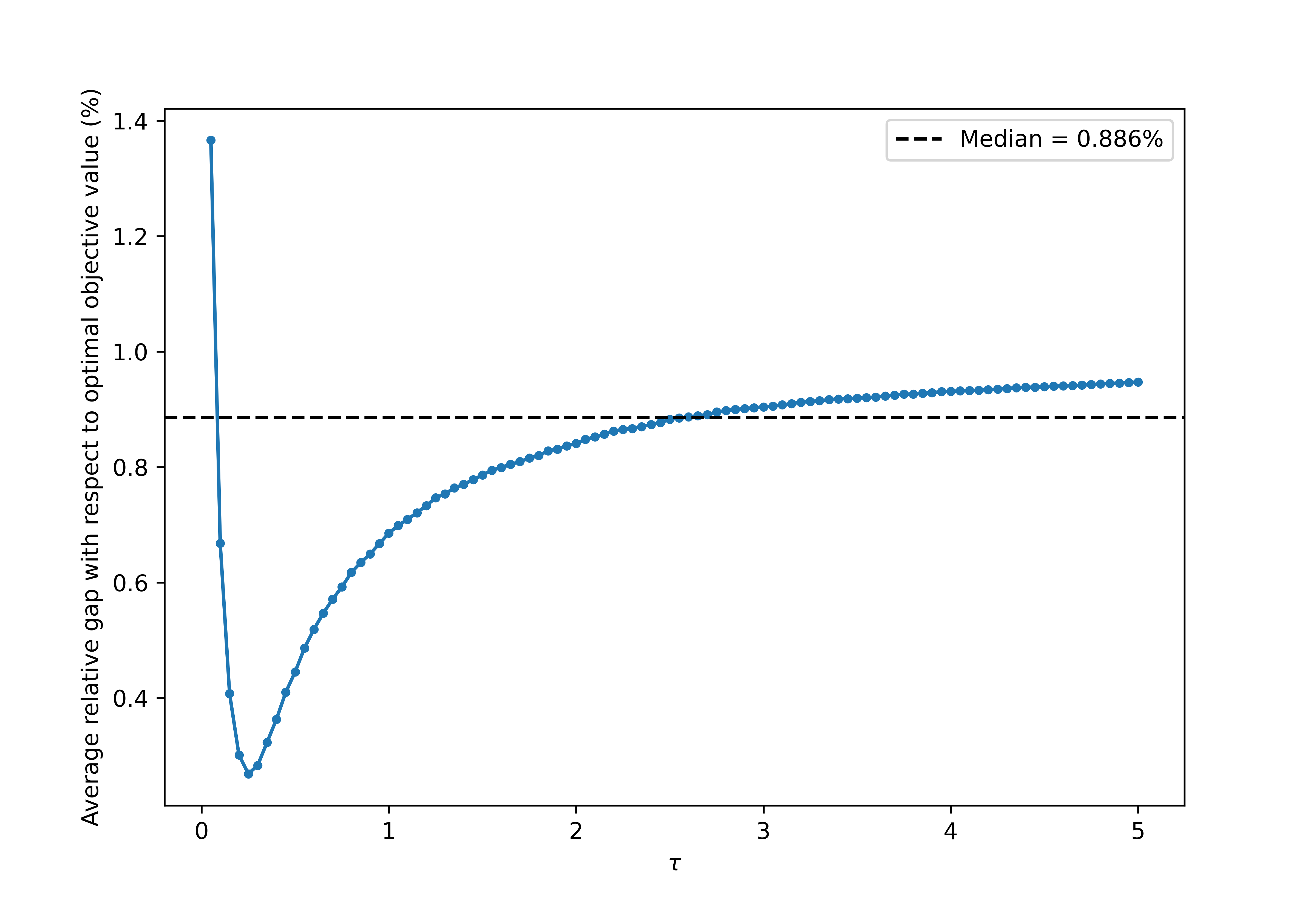}
        \label{fig:tauC_2}
    }
    \caption{Distribution of optimal $\tau$ in terms of tightness and solution optimality.}
\end{figure}

\hyperref[fig:tauC_1]{Figure \ref*{fig:tauC_1}} presents the distribution of the optimal $\tau$ values, where the black dashed vertical line indicates the median of optimal $\tau$. We observe that, in most cases, a small value of $\tau$ tends to yield the tightest approximation. Specifically, the optimal value of $\tau$ is always smaller than 0.05 and falls below 0.02 for more than 50\% of the cases.  Because the tightness-optimal $\tau$ is not necessarily optimal for the true solution performance, we also search for the optimal $\tau\in\{0.05,0.10,\ldots,5.00\}$ such that the resulting solution attains the highest expected utility.\footnote{We start from 0.05 because extremely small values $\tau$ may lead to numerical instability in our later optimization.} For each instance, we calculate the relative gap between the objective value associated with a given $\tau$ and that achieved under the optimal choice $\tau^{\star}$, which assesses the optimality gap for a specific $\tau$. \hyperref[fig:tauC_2]{Figure \ref*{fig:tauC_2}} reports the average optimality gap of different values of $\tau$, where the horizontal dashed line indicates the median of the average gap. It shows that, as $\tau$ increases, the optimality gap first drops sharply and then rises gradually. Overall, the best-performing value of $\tau$ that leads to the best matching solution is around $0.25$, and the optimality gap remains small for a range of values of $\tau$. Specifically, except for $\tau=0.05$, the average optimality gap is smaller than 1\% when $\tau$ is not so large.

In summary, we find that a small value of $\tau$ (but not too small to cause numerical issues) is generally favorable in terms of both approximation and optimality gap. In such a regime, the first (exponential) term in the bound function of \hyperref[corollary1]{Corollary \ref*{corollary1}} dominates and diminishes the role of the other terms. This motivates us to focus on a simple approximation scheme centered around the exponential function. It turns out that this approximation enjoys strong theoretical guarantees, which will be discussed in detail in the following sections.

\subsubsection{Approximation Formulation.}
Motivated by the above upper bound and our numerical findings, we propose the following approximation model (A-SP) as a proxy optimization problem for SP-SAA:
\begin{subequations}
    \begin{align}
       \textbf{[A-SP]}:  \max \ & \tau\sum_{i \in \mathcal{A}^{D}}  \log \sum_{j \in \mathcal{A}^{S}} \exp{\left(\frac{u_{ij}}{\tau}\right)} p_{ij}x_{ij} \label{eq:5a}\\
       {\rm s.t.} \ & \hyperref[eq:a1]{(\ref*{eq:a1})-(\ref*{eq:a3})} \nonumber,
    \end{align}
\end{subequations}

\noindent where the objective function \hyperref[eq:5a]{(\ref*{eq:5a})} remains concave and is close to the original upper bound for a small value of $\tau$, which is the regime where the approximation works well. We observe that A-SP is a mixed-integer convex optimization problem, which can be solved by outer approximation or the branch-and-bound algorithm. To improve the efficiency of the approximation model, in the following proposition, we show that A-SP admits an equivalent mixed-integer exponential cone program (MIECP) reformulation that can be solved by specialized solvers such as Mosek.

\begin{Proposition}\label{pro:miecp}
The optimal value and the set of optimal solutions of model A-SP coincide with those of the following mixed-integer exponential cone program:
\begin{subequations}
    \begin{align}
        \textbf{[ECA-SP]}: \max \ &\tau \sum_{i \in \mathcal{A}^{D}}Z_{i} \\
        {\rm s.t.} \ & \hyperref[eq:a1]{(\ref*{eq:a1})-(\ref*{eq:a3})}, \nonumber \\
        & Y_{i} = \sum_{j \in \mathcal{A}^{S} }\exp{\left(\frac{u_{ij}}{\tau}\right)}p_{ij}x_{ij}, \forall i \in \mathcal{A}^{D},\\
        & \left( Y_{i},1,Z_{i}\right) \in \mathcal{K}_{\exp},Y_{i}\ge 0, \forall i \in \mathcal{A}^{D}.
    \end{align}
\end{subequations}
\end{Proposition}


\begin{Remark}
When the supply is very limited, there may exist a demand $i^{\prime}$ such that $x_{i^{\prime}j}=0$ for all $j \in \mathcal{A}^{S}$, which would make $Y_{i^{\prime}}=0$. In such cases, one can introduce a small enough constant ($\epsilon>0$), so that $Y_{i^{\prime}}= \sum_{j \in \mathcal{A}^{S} }\exp{\left(\frac{u_{i^{\prime}j}}{\tau}\right)}p_{i^{\prime}j}x_{i^{\prime}j}+\epsilon$ and the function is well defined everywhere. Our numerical test suggests that the addition of $\epsilon$ would not affect the optimality of the solution. 
\end{Remark}

\subsubsection{Approximation Performance Analysis.} \label{sec:analysis}
Next, we evaluate the approximation performance of the optimal solution of ECA-SP, which is denoted as $\mathbf{x}^{A}$. Accordingly, the approximation gap of ECA-SP in terms of solving the original stochastic program SP is computed as $\frac{Obj^{E} (\mathbf{x}^{\star})-Obj^{E}(\mathbf{x}^{A})}{Obj^{E} (\mathbf{x}^{\star})}$ (note that the objective function of ECA-SP is denoted as $Obj^{A}(\mathbf{x})$ instead). Let $\underline{p} = \min_{i,j} p_{ij}$ and $\overline{p} = \max_{i,j} p_{ij}$. 
We first introduce the following lemma. 
\begin{Lemma}
\label{lemma:bounds}
Given any feasible solution $\mathbf{x}$, let $V_i^{E}(\mathbf{x}) = \mathbb{E}\left[\max_{j\in \mathcal{A}^S} u_{ij}\tilde{\xi}_{ij}x_{ij}\right]$ and $V_{i}^{A}(\mathbf{x}) = \tau \log \left( \sum_{j \in \mathcal{A}^{S}} p_{ij} \exp{\left(\frac{u_{ij}}{\tau}\right)}x_{ij} \right)$ so that $Obj^{E}(\mathbf{x})=\sum_{i \in \mathcal{A}^{D}} V_{i}^{E}(\mathbf{x})$ and $Obj^{A}(\mathbf{x})=\sum_{i \in \mathcal{A}^{D}} V_{i}^{A}(\mathbf{x})$. Then the following inequalities hold:
\begin{align*}
    & V_{i}^{E}(\mathbf{x}) \ge \underline{p}U_{i}(\mathbf{x}),\\
    &  V_{i}^{E}(\mathbf{x}) \le \overline{q}U_{i}(\mathbf{x}), \\
    &  V_{i}^{A}(\mathbf{x}) \ge \tau \log(\underline{p})+ U_{i}(\mathbf{x}),   \\
    & V_{i}^{A}(\mathbf{x})  \le \tau \log(\overline{p})+ U_{i}(\mathbf{x})+ \tau \log (\theta)  ,
\end{align*}
\noindent where $\overline{q}=\left(1-(1-\overline{p})^{\theta}\right)$ and  $U_{i}(\mathbf{x})=\max_{j\in \mathcal{A}^S}u_{ij}x_{ij}$ for  $i\in \mathcal{A}^{D}$.
\end{Lemma}

Based on \hyperref[lemma:bounds]{Lemma \ref*{lemma:bounds}}, we first obtain the performance guarantee of ECA-SP when the utility values are drawn from uniform distributions, the supply-demand ratio is equal to $\theta$, and the acceptance probabilities are homogeneous. Note that the uniform distribution assumption is mainly technical to enable a cleaner expression of the bound and will be relaxed later. 


\begin{Theorem}\label{theoemhomo}
Suppose the utility values $u_{ij}$ are i.i.d. random variables following a uniform distribution in $[a_{i},b_{i}]$ for $i\in \mathcal{A}^D$ and $\gamma = \theta$. When $p_{ij}=p$ for all $i\in\mathcal{A}^D$ and $j\in\mathcal{A}^S$, the expected approximation gap of ECA-SP satisfies
    \begin{align*}
        \frac{\mathbb{E}\left[Obj^{E} (\mathbf{x}^{\star})-Obj^{E}(\mathbf{x}^{A})\right]}{\mathbb{E}\left[Obj^{E} (\mathbf{x}^{\star})\right]} \leq \frac{\tau |\mathcal{A}^{D}| \log (\theta)+\left(1-\frac{p}{q} \right)\sum_{i \in \mathcal{A}^{D}}\left(b_i-a_i\right)}{\sum_{i \in \mathcal{A}^{D}}a_i+\left(\frac{\frac{\theta}{q}-\frac{1}{p}+1}{\theta+1}\right)\sum_{i \in \mathcal{A}^{D}}\left(b_i-a_i\right)}
    \end{align*}
    \noindent where $q=\left(1-\left(1-p \right) ^{\theta}\right)$.
\end{Theorem}

\hyperref[theoemhomo]{Theorem \ref*{theoemhomo}} has several important implications. First, the gap is logarithmic in $\theta$, so the bound would be tighter for smaller $\theta$ (as shown below, the gap is zero if $\theta=1$); however, one can choose a small $\tau$ to offset the impact of $\theta$. Recall that a small $\tau$ is preferred as we discussed in Section \ref{sec:choice}. In fact, as $\tau\rightarrow 0$, the gap is approaching zero for $p=1$, which corresponds to the case without acceptance uncertainty. Second, as one might expect, the bound is tighter as the utility ranges become narrower and the acceptance probability increases. Specifically, consider the case where $p$ is close to one and $\theta$ is large, then the approximation gap is at most
\begin{equation*}
    \frac{(1-\frac{p}{q})(\theta+1)}{\frac{\theta}{q}-\frac{1}{p}+1} \approx \frac{(1-p)(\theta+1)}{\theta - \frac{1}{p}+1} \approx 1-p,
\end{equation*}
which decreases linearly in the acceptance probability. This parametric bound does not approach one as $\theta$ becomes large, which stands in contrast to the limitation of DAP in the big-$\theta$ regime. For a practical case where $\theta=4$,  $a=5$, $b=10$, $p=0.8$ and $\tau=0.1$, the derived worst-case bound is 12.93\% (the bound is not necessarily tight, and numerically, we find the gap to be much smaller, as shown later in Section \ref{sec:results}). 
 The following result confirms that when $\theta=1$, the proposed MIECP is exact. 

\begin{Corollary}\label{cor1} When the acceptance probabilities are homogeneous and the maximum number of recommendations $\theta$ equals 1, ECA-SP is exact.
\end{Corollary}

Additionally, in the following, we provide a sufficient condition under which the ECP-based policy dominates the DIMP policy. Note that when $\theta = 1$, DIMP is also exact and attains the optimal solution of the original stochastic program.

\begin{Proposition}\label{dominate}
Based on the conditions in \hyperref[theoemhomo]{Theorem \ref*{theoemhomo}}, we further assume that $b_i = (1+\alpha)a_i$. When $\tau \rightarrow 0$, for $\theta \ge 2$, there exists $\alpha_{\max}(\theta,p)=\frac{(q-p)(\theta+1)}{\,q/p-\theta(1-p)^\theta}>0$ such that for any $0<\alpha<\alpha_{\max}(\theta,p)$, the MIECP solution is superior to the DIMP solution, i.e., $\mathbb{E}\left[Obj^{E}(\mathbf{x}^{A})-Obj^{E} (\mathbf{x}^{M})\right]\ge 0$.
\end{Proposition}

\hyperref[dominate]{Proposition \ref*{dominate}} implies that for most relevant cases ($\theta\geq2$), one can expect the ECP-based matching solution to prevail against the LP-based approximation. The condition on $\alpha$ is relatively mild: for example, when $\theta=4$ and $p=0.8$, we have $\alpha_{\max}(\theta,p)=0.8$, and the utilities can be drawn from a wide range in $[a_i, 1.8a_i]$ for the dominance to hold. This extends the worst-case result to stronger average-case comparisons.

Next, we analyze the approximation performance under general settings where the acceptance probabilities are heterogeneous and the utility values are generally distributed. 
Assume that the utility of demand $i$ has a lower bound $a_i$ for any supplier, the following bound holds when we still assume independence of supplier acceptances facing the same demand (while there could be correlation of supplier acceptance across demands).

\begin{Theorem}\label{theoremheter}
Suppose there exists $a_i$ such that $u_{ij}\geq a_i$ for $i\in\mathcal{A}^D$, the approximation gap of ECA-SP satisfies
$$\frac{Obj^{E} (\mathbf{x}^{\star})-Obj^{E} (\mathbf{x}^{A})}{Obj^{E} (\mathbf{x}^{\star})}\le   1-\frac{\underline{p}}{\overline{q}}+\frac{\tau\underline{p}\left| \mathcal{A}^{D}\right| \log(\frac{\theta\overline{p}}{\underline{p}})}{\sum_{i\in \mathcal{A}^{D\prime}}\left(1-(1-\underline{p})^{\lfloor \gamma \rfloor+1}\right)a_{i}+ \sum_{i\in \mathcal{A}^{D}\setminus \mathcal{A}^{D\prime}}\left(1-(1-\underline{p})^{\min(\lfloor \gamma \rfloor, \theta)}\right)a_{i}}$$
\noindent where $\mathcal{A}^{D\prime}$ is the set of $\left(|\mathcal{A}^{S}|-\lfloor \gamma \rfloor|\mathcal{A}^{D}|\right)^{+}$ demands with the highest lower bound of the utility when $\gamma<\theta$, otherwise, it is the empty set; and $\overline{q}=\left(1-(1-\overline{p})^{\theta}\right)$.
\end{Theorem}

Compared to \hyperref[theoemhomo]{Theorem \ref*{theoemhomo}}, the general performance bound derived in \hyperref[theoremheter]{Theorem \ref*{theoremheter}} is less tight but shares a similar structure that also includes a logarithmic term in $\theta$. The logarithmic term can be controlled by $\tau$ so that the approximation gap is determined mainly by $1-\underline{p}/\overline{q}$, which decreases in $\underline{p}$ and increases in $\bar{p}$. Similar to  \hyperref[theoemhomo]{Theorem \ref*{theoemhomo}}, as $\tau\rightarrow0$, the approximation gap reduces to zero when $\underline{p}=\overline{p}=1$. In a practical case where $\gamma=\theta=4$,  $a_i=5$, $p\in [0.7, 0.9]$ and $\tau=0.1$, the derived worst-case bound is 32.30\%. 

Furthermore, following \hyperref[theoremheter]{Theorem \ref*{theoremheter}}, we can prove a similar bound even if there is a correlation in the suppliers' acceptance facing the same demand. Intuitively, when correlation arises, the simultaneous rejection of demands from suppliers is more likely, increasing the risk of failed matches. 
\begin{Corollary}\label{corogeneral}
    When the supplier's response is correlated, the approximation gap of ECA-SP satisfies
    \[\frac{Obj^{E} (\mathbf{x}^{\star})-Obj^{E} (\mathbf{x}^{A})}{Obj^{E} (\mathbf{x}^{\star})}\le   1-\underline{p}+\frac{\tau\left| \mathcal{A}^{D}\right| \log(\frac{\theta\overline{p}}{\underline{p}})}{\sum_{i\in \mathcal{A}^{\underline{D}}}a_{i}}, \]
    where $\mathcal{A}^{\underline{D}}=\mathcal{A}^{D\prime}$ when $\lfloor \gamma \rfloor=0$; otherwise, $\mathcal{A}^{\underline{D}}=\mathcal{A}^{D}$.
\end{Corollary}

We note that the major term of the above bound is linear in $\underline{p}$, and $\bar{p}$ plays a relatively minor role. Overall, the general bound presented in \hyperref[corogeneral]{Corollary \ref*{corogeneral}}, albeit less tight, follows the structure of the bounds under more specialized scenarios. It suggests the ECP-based policy's versatility. Moreover, all the performance analysis indicates that a small $\tau$ is preferred to reduce the worst-case bound, which corroborates with our findings in Section \ref{sec:choice}. However, if $\tau$ is too small, the exponential term in the objective function of MIECP will grow too large to admit reliable solutions in numerical solvers. Therefore, for practical use, we suggest a modestly small value of $\tau$, such as 0.1, to balance numerical stability and performance.

\section{Numerical Experiments}\label{sec:results}
In this section, we evaluate the performance of the proposed approximation algorithm against benchmarks on both synthetic and real-world data sets. The algorithms were implemented in Python using Gurobi 12 for MILP and Mosek 11 for MIECP. All numerical experiments were conducted on a Windows 10 machine with a 64-bit operating system and 16GB of RAM. The computation time is measured in seconds, with a time limit set to 3,600 seconds.

\subsection{Synthetic Experiment} \label{sec:comp}
We generate test instances of varying sizes. As shown in \hyperref[tab:homo]{Table \ref*{tab:homo}}, instance D$i$-S$j$ includes $i$ demands and $j$ suppliers. The platform utility is assumed to take the form of 
$u_{ij}= 0.4 + 0.2u_i^{D} + 0.2u_j^{S} + 0.2u_{ij}^{R}$, where the three terms denote demand, supply, and pair-specific utility values, and they are drawn uniformly within $[0,1]$. For each problem size, we create 10 random instances by drawing the utility values from the specified distribution. Besides, we let $\theta=4$ and $\tau=0.1$ following our earlier discussion on the choice of $\tau$.

\subsubsection{Homogeneous acceptance probabilities.} \hyperref[tab:homo]{Table \ref*{tab:homo}} summarizes the evaluation results for problems with homogeneous acceptance probabilities ($p=0.8$). Rows `Gap-A (\%)' and `Gap-W (\%)' represent the average and maximum percentage differences, respectively, in the objective value between each formulation and the best objective value found across the four approaches. Rows `CPU' report the average computation time (in seconds) across the ten runs. Note that because reformulation R-SP is exact, its gap is zero unless the problem size is too large to yield exact optimal solutions (e.g., in the last row of \hyperref[tab:homo]{Table \ref*{tab:homo}}). Particularly, R-SP can be computationally expensive for large instances, as the running time often exceeds the set time limit. In contrast, MIECP consistently delivers high-quality solutions (with an average gap of less than 2\%) using only a small fraction of the computational time of R-SP. We also observe that although DAP and DIMP request the least computation time, their optimality gaps can be huge. For DAP, the gap exceeds 30\% in several cases. This loss tends to be larger when $\gamma$ is smaller, whereas the performance improves as $\gamma$ increases, which is consistent with \hyperref[prop:unbound]{Proposition \ref*{prop:unbound}}. For DIMP, except for the cases with $\gamma=1$, the average gap is typically above 15\%.
Overall, although R-SP is exact with homogeneous acceptance probabilities, MIECP can be preferable to achieve a short matching time when handling large instances. 


\begin{table}[htbp]
    \caption{Performance evaluation with homogeneous acceptance probabilities }
    \label{tab:homo}
    \centering
    \small
    {\textbf{(a) Small Instances} \par}\vspace{0.3em}
    \begin{tabularx}{0.9\textwidth}{cc RRRRR}
    \hline
    \multirow{2}{*}{Metric} & \multirow{2}{*}{Method} & \multicolumn{5}{c}{Instance} \\ \cline{3-7}
    & & D10-S10 & D10-S20 & D10-S30 & D10-S40 & D10-S50 \\ \hline
    \multirow{4}{*}{Gap-A (\%)} 
    & DAP   & 43.94 & 29.85 & 13.63 & 0.85 & 0.70 \\
    & DIMP  & 0.00  & 14.96 & 17.89 & 18.84 & 18.95 \\
    & R-SP  & 0.00  & 0.00  & 0.00  & 0.00  & 0.00 \\
    & MIECP & 0.00  & 0.88  & 1.17  & 0.62  & 0.18 \\ \hline
    \multirow{4}{*}{Gap-W (\%)} 
    & DAP   & 51.22 & 38.43 & 18.33 & 1.22 & 1.13 \\
    & DIMP  & 0.00  & 15.23 & 18.21 & 19.00 & 19.08 \\
    & R-SP  & 0.00  & 0.00  & 0.00  & 0.00  & 0.00 \\
    & MIECP & 0.00  & 1.80  & 2.60  & 1.12  & 0.45 \\ \hline
    \multirow{4}{*}{CPU}        
    & DAP   & 0.01  & 0.01  & 0.01  & 0.02  & 0.03 \\
    & DIMP  & 0.00  & 0.01  & 0.01  & 0.01  & 0.01 \\
    & R-SP  & 0.05  & 6.58  & 1300.65 & 3600.00 & 3600.00 \\
    & MIECP & 0.03  & 0.04  & 0.07  & 0.09  & 0.43 \\ \hline
    \end{tabularx}

    \vspace{1em}
    
    {\textbf{(b) Moderate Instances} \par}\vspace{0.3em}
    \begin{tabularx}{0.9\textwidth}{cc RRRRR}
    \hline
    \multirow{2}{*}{Metric} & \multirow{2}{*}{Method} & \multicolumn{5}{c}{Instance} \\ \cline{3-7}
    & & D20-S50 & D20-S80 & D30-S60 & D30-S75 & D30-S100 \\ \hline
    \multirow{4}{*}{Gap-A (\%)} 
    & DAP   & 25.99 & 1.11  & 39.12 & 28.15 & 11.88 \\
    & DIMP  & 16.63 & 18.67 & 15.04 & 16.63 & 18.22 \\
    & R-SP  & 0.00  & 0.00  & 0.00  & 0.00  & 0.00 \\
    & MIECP & 0.79  & 0.26  & 3.37  & 1.46  & 0.72 \\ \hline
    \multirow{4}{*}{Gap-W (\%)} 
    & DAP   & 30.55 & 1.51  & 44.05 & 31.54 & 15.21 \\
    & DIMP  & 16.77 & 18.88 & 15.21 & 16.81 & 18.38 \\
    & R-SP  & 0.00  & 0.00  & 0.00  & 0.00  & 0.00 \\
    & MIECP & 1.93  & 0.51  & 4.26  & 2.38  & 1.01 \\ \hline
    \multirow{4}{*}{CPU}        
    & DAP   & 0.03  & 0.04  & 0.04  & 0.04  & 0.05 \\
    & DIMP  & 0.01  & 0.02  & 0.02  & 0.02  & 0.03 \\
    & R-SP  & 3600.00 & 3600.00 & 3600.00 & 3600.00 & 3600.00 \\
    & MIECP & 3.91  & 4.92  & 33.84 & 15.53 & 7.38 \\ \hline
    \end{tabularx}

    \vspace{1em}
    
    {\textbf{(c) Large Instances} \par}\vspace{0.3em}
    \begin{tabularx}{0.9\textwidth}{cc RRRRR}
    \hline
    \multirow{2}{*}{Metric} & \multirow{2}{*}{Method} & \multicolumn{5}{c}{Instance} \\ \cline{3-7}
    & & D50-S150 & D50-S200 & D50-S250 & D100-S250 & D100-S400 \\ \hline
    \multirow{4}{*}{Gap-A (\%)} 
    & DAP   & 20.33 & 1.22  & 1.00  & 31.11 & 1.10 \\
    & DIMP  & 17.93 & 18.65 & 18.94 & 16.42 & 18.48 \\
    & R-SP  & 0.00  & 0.00  & 0.00  & 0.00  & 0.05 \\
    & MIECP & 1.36  & 0.34  & 0.34  & 1.97  & 0.12 \\ \hline
    \multirow{4}{*}{Gap-W (\%)} 
    & DAP   & 21.65 & 1.63  & 1.20  & 31.87 & 1.53 \\
    & DIMP  & 18.06 & 18.77 & 19.01 & 16.57 & 18.58 \\
    & R-SP  & 0.00  & 0.00  & 0.00  & 0.00  & 0.45 \\
    & MIECP & 1.80  & 0.46  & 0.43  & 2.66  & 0.28 \\ \hline
    \multirow{4}{*}{CPU}        
    & DAP   & 0.10  & 0.12  & 0.15  & 0.34  & 0.58 \\
    & DIMP  & 0.07  & 0.07  & 0.09  & 0.16  & 0.25 \\
    & R-SP  & 3600.00 & 3548.05 & 3600.00 & 3600.00 & 3600.00 \\
    & MIECP & 8.11  & 5.45  & 3.01  & 82.97 & 8.36 \\ \hline
    \end{tabularx}
    
\end{table}

\begin{figure}[htbp]
    \centering
    \subfigure[Performance comparison with varying $\theta$]{
        \includegraphics[width=0.48\textwidth]{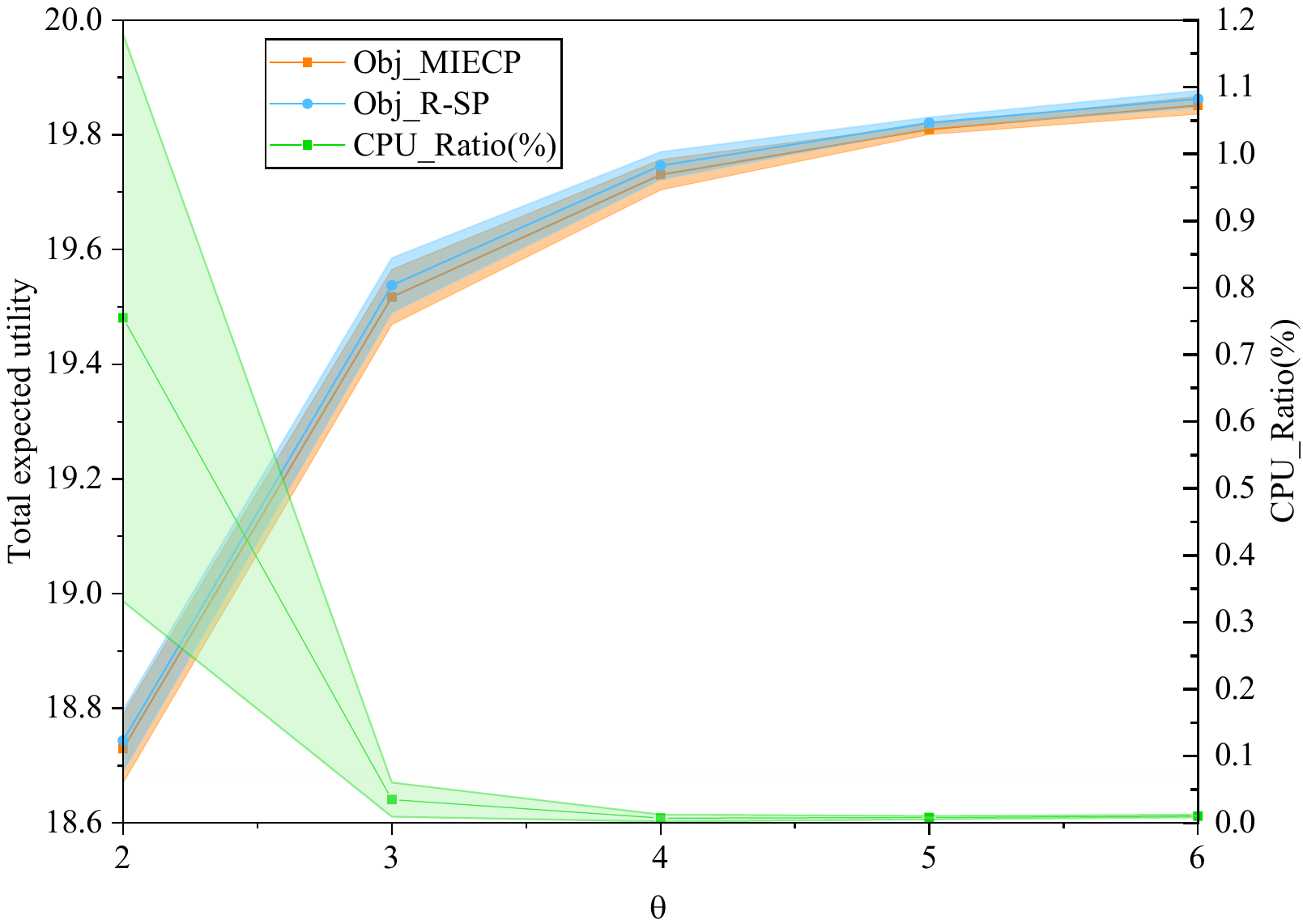}
        \label{fig:subsent}
    }
    \hfill
    \subfigure[Performance comparison with varying $p$]{
        \includegraphics[width=0.48\textwidth]{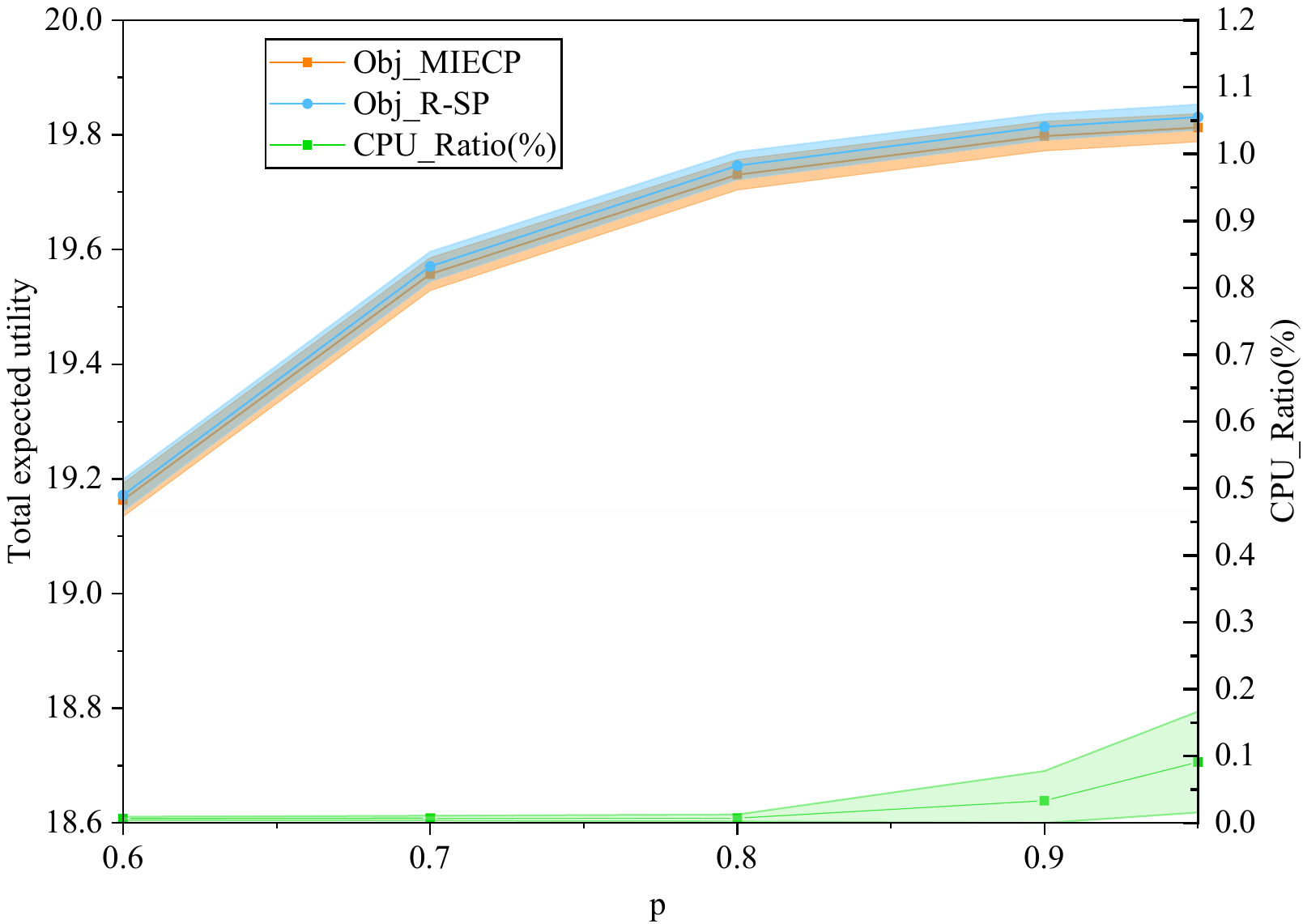}
        \label{fig:subsenp}
    }
    \caption{Optimal objective values and running time of R-SP and MIECP on D20-S80 instances}
    \label{fig:sensitive}
\end{figure}

To understand the impact of the matching conditions, we fix $\gamma=\theta$ and evaluate the performance of R-SP and MIECP under different values of $\theta$ and $p$, while the utility values are randomly drawn from $[0.4,1]$. \hyperref[fig:sensitive]{Figure \ref*{fig:sensitive}} presents the optimal objective values and running time of R-SP and MIECP (CPU\_Ratio (\%) indicates the ratio of the running time of MIECP to that of R-SP), including the mean values and 95\% confidence intervals in shaded areas. We observe that MIECP maintains near optimality across different parameter regimes (within 0.1\%), whereas its computation time can be one to two orders shorter than R-SP. Notably, when $\theta$ is large or $p$ is small (the more challenging cases), the required running time of MIECP is less than 0.1\% that of R-SP.

\subsubsection{Heterogeneous acceptance probabilities.} Next, for instances with heterogeneous acceptance probabilities, in addition to DAP, DIMP, SAA, and MIECP, we consider an extension of the matching policy in \citet{ekbatani2026lyft1} for a notification problem that shares a similar structure to ours but without capacity constraints. Specifically, we extend the submodular welfare maximization algorithm to incorporate the recommendation cap, and the resulting benchmark is referred to as Cap-SW. 
In particular, to respect the recommendation cap, we need to apply a continuous-greedy type fractional approximation over the corresponding down-monotone polytope, followed by a new lossless rounding procedure along simple paths/cycles on the bipartite graph. Due to the limited discretization accuracy, the implementable guarantee takes the standard form \((1-1/e-\varepsilon)\) for any prescribed \(\varepsilon>0\). The algorithm pseudocode is presented in \ref{sec:sub}, and \hyperref[propsub]{Proposition \ref*{propsub}} formalizes the performance guarantee of Cap-SW.
\begin{Proposition}\label{propsub}
    For any prescribed \(\varepsilon>0\), consider the problem in which each supplier can be recommended to at most one demand and each demand can be recommended to at most $\theta$ suppliers. 
    Then there exists a randomized polynomial-time algorithm that returns a feasible recommendation solution $\mathbf{x}^{S}$ such that
\[
\mathbb E\left[Obj^{E}(\mathbf{x}^{S})\right]\ge
\left(1-1/e-\varepsilon\right)Obj^{E} (\mathbf{x}^{\star}).
\]
\end{Proposition}


We also note that one may apply an L-shaped decomposition method to more effectively solve SP-SAA. In \ref{addLshaped}, we implement and test the L-shaped method, where we further tighten the formulation with additional acceleration techniques. However, we did not observe significant computational benefits from the L-shaped method over Gurobi's default MILP solver. Therefore, our SAA implementation evaluation result is still based on Gurobi.  

The results for the instances with heterogeneous acceptance probabilities are shown in \hyperref[tab:heter]{Table \ref*{tab:heter}}, where the acceptance probability $p_{ij}$ is drawn from $[0.7,0.9]$. We set the number of scenarios in the SAA method to 1,000 for small-scale instances and 100 for moderate- and large-scale instances (a larger number would exceed the time limit). Note that SAA may not always find the optimal solution even for small instances, and we observe an average optimality gap of 0.08\%, compared with 0.80\%, 14.45\%, 11.12\%, and 3.63\% for MIECP, DAP, DIMP, and Cap-SW, respectively. For moderate-to-large problems, the MIECP approximation achieves comparable objective values to those obtained by SAA, while being more than 60 times faster than SAA. Particularly, MIECP could outperform SAA in both solution quality and computation time on large instances. 

\begin{figure}[htbp]
    \centering
    \subfigure[Performance comparison with varying $\theta$]{
        \includegraphics[width=0.48\textwidth]{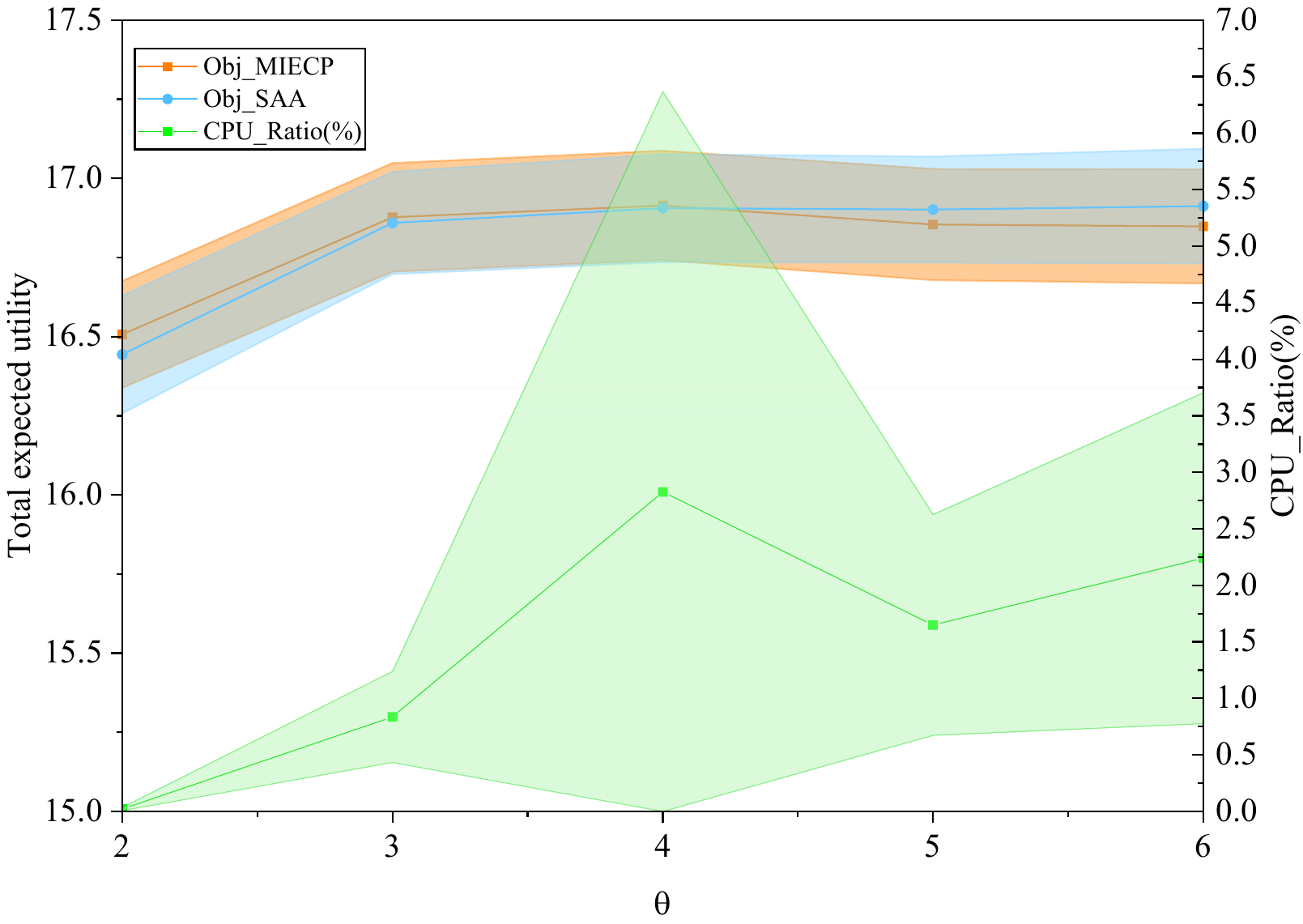}
        \label{fig:subrobusttheta}
    }
    \hfill
    \subfigure[Performance comparison with varying $\gamma$]{
        \includegraphics[width=0.48\textwidth]{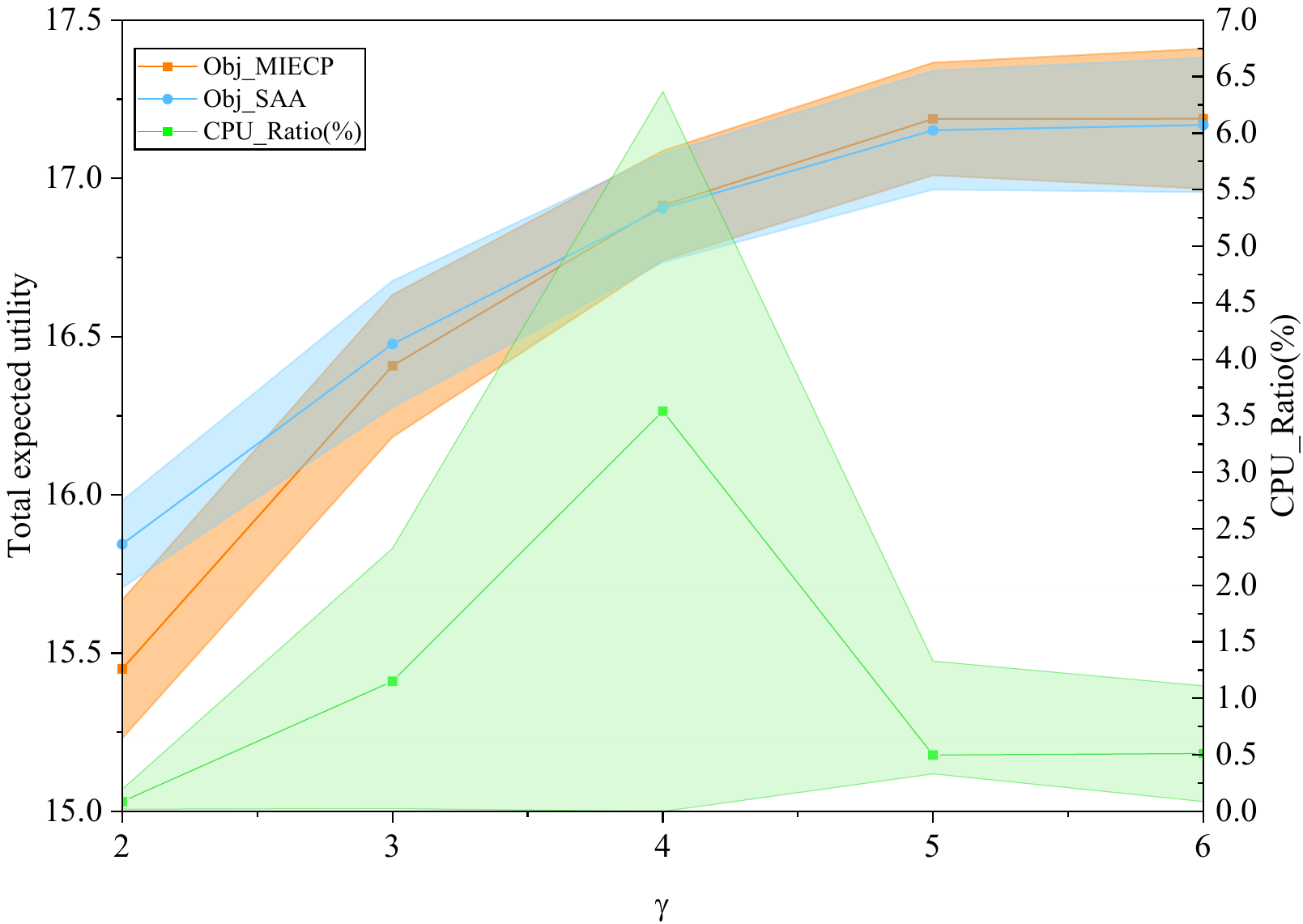}
        \label{fig:subrobustgamma}
    }
    \caption{Optimal objective values and running time of SAA and MIECP with $|\mathcal{A}^{D}|=20$}
    \label{fig:robust}
\end{figure}

\begin{table}[htbp]
    \caption{Performance evaluation with heterogeneous acceptance probabilities}
    \label{tab:heter}
    \centering
    \small
    
    { \textbf{(a) Small Instances} \par}\vspace{0.3em}
    \begin{tabularx}{0.9\textwidth}{cc RRRRR}
    \hline
    \multirow{2}{*}{Metric} & \multirow{2}{*}{Method} & \multicolumn{5}{c}{Instance} \\ \cline{3-7}
    & & D10-S10 & D10-S20 & D10-S30 & D10-S40 & D10-S50 \\ \hline
    \multirow{5}{*}{Gap-A (\%)} 
    & DAP    & 41.50 & 26.58 & 9.52  & 2.50  & 1.62 \\
    & DIMP   & 0.00  & 11.62 & 13.45 & 14.29 & 14.14 \\
    & SAA    & 0.26  & 0.00  & 0.00  & 0.00  & 0.00 \\
    & Cap-SW & 0.88  & 1.93  & 2.24  & 1.60  & 1.01 \\
    & MIECP  & 0.14  & 1.99  & 2.09  & 2.58  & 0.12 \\ \hline
    \multirow{5}{*}{Gap-W (\%)} 
    & DAP    & 50.80 & 36.74 & 17.89 & 3.55  & 2.90 \\
    & DIMP   & 0.00  & 12.53 & 14.49 & 15.12 & 15.97 \\
    & SAA    & 0.53  & 0.00  & 0.00  & 0.00  & 0.02 \\
    & Cap-SW & 2.52  & 3.65  & 3.96  & 2.61  & 1.97 \\
    & MIECP  & 0.83  & 3.77  & 4.12  & 3.63  & 0.23 \\ \hline
    \multirow{5}{*}{CPU}        
    & DAP    & 0.00  & 0.01  & 0.01  & 0.01  & 0.01 \\
    & DIMP   & 0.01  & 0.00  & 0.00  & 0.01  & 0.01 \\
    & SAA    & 3.35  & 2956.77 & 3283.65 & 488.24 & 491.18 \\
    & Cap-SW & 3.72  & 5.02  & 5.84  & 7.04  & 8.06 \\
    & MIECP  & 0.03  & 0.06  & 0.09  & 0.06  & 0.06 \\ \hline
    \end{tabularx}

    \vspace{1em}
    
    {\textbf{(b) Moderate Instances} \par}\vspace{0.3em}
    \begin{tabularx}{0.9\textwidth}{cc RRRRR}
    \hline
    \multirow{2}{*}{Metric} & \multirow{2}{*}{Method} & \multicolumn{5}{c}{Instance} \\ \cline{3-7}
    & & D20-S50 & D20-S80 & D30-S60 & D30-S75 & D30-S100 \\ \hline
    \multirow{5}{*}{Gap-A (\%)} 
    & DAP    & 21.51 & 1.98  & 32.51 & 23.20 & 7.47 \\
    & DIMP   & 11.79 & 12.98 & 9.66  & 11.15 & 12.01 \\
    & SAA    & 0.00  & 0.11  & 0.00  & 0.01  & 0.00 \\
    & Cap-SW & 4.00  & 2.31  & 4.59  & 4.83  & 4.09 \\
    & MIECP  & 1.00  & 0.06  & 1.81  & 0.86  & 0.20 \\ \hline
    \multirow{5}{*}{Gap-W (\%)} 
    & DAP    & 27.88 & 2.66  & 39.03 & 28.60 & 11.44 \\
    & DIMP   & 12.33 & 14.19 & 10.27 & 11.72 & 12.92 \\
    & SAA    & 0.01  & 0.43  & 0.00  & 0.08  & 0.04 \\
    & Cap-SW & 5.09  & 3.39  & 6.07  & 7.02  & 5.20 \\
    & MIECP  & 1.73  & 0.46  & 2.46  & 1.83  & 0.48 \\ \hline
    \multirow{5}{*}{CPU}        
    & DAP    & 0.02  & 0.02  & 0.03  & 0.03  & 0.04 \\
    & DIMP   & 0.01  & 0.01  & 0.01  & 0.04  & 0.10 \\
    & SAA    & 1866.12 & 68.90 & 3600.00 & 3508.59 & 1613.77 \\
    & Cap-SW & 16.59 & 22.98 & 28.42 & 32.74 & 39.83 \\
    & MIECP  & 1.92  & 1.05  & 14.58 & 11.53 & 2.32 \\ \hline
    \end{tabularx}

    \vspace{1em}
    
    { \textbf{(c) Large Instances} \par}\vspace{0.3em}
    \begin{tabularx}{0.9\textwidth}{cc RRRRR}
    \hline
    \multirow{2}{*}{Metric} & \multirow{2}{*}{Method} & \multicolumn{5}{c}{Instance} \\ \cline{3-7}
    & & D50-S150 & D50-S200 & D50-S250 & D100-S250 & D100-S400 \\ \hline
    \multirow{5}{*}{Gap-A (\%)} 
    & DAP    & 16.25 & 1.73  & 1.72  & 27.02 & 1.61 \\
    & DIMP   & 11.36 & 11.77 & 11.96 & 9.51  & 11.07 \\
    & SAA    & 0.00  & 0.28  & 0.28  & 0.02  & 0.31 \\
    & Cap-SW & 5.95  & 4.26  & 2.45  & 8.26  & 6.09 \\
    & MIECP  & 0.52  & 0.00  & 0.00  & 0.58  & 0.00 \\ \hline
    \multirow{5}{*}{Gap-W (\%)} 
    & DAP    & 19.51 & 2.32  & 2.13  & 29.50 & 2.02 \\
    & DIMP   & 11.78 & 12.28 & 12.16 & 10.10 & 11.26 \\
    & SAA    & 0.00  & 0.51  & 0.49  & 0.24  & 0.47 \\
    & Cap-SW & 7.92  & 5.17  & 3.07  & 9.86  & 7.38 \\
    & MIECP  & 0.78  & 0.00  & 0.00  & 1.33  & 0.00 \\ \hline
    \multirow{5}{*}{CPU}        
    & DAP    & 0.07  & 0.10  & 0.15  & 0.36  & 0.64 \\
    & DIMP   & 0.18  & 0.27  & 0.34  & 0.63  & 0.69 \\
    & SAA    & 3600.00 & 3600.00 & 3600.00 & 3600.00 & 3600.00 \\
    & Cap-SW & 90.84 & 116.65 & 138.56 & 281.16 & 455.22 \\
    & MIECP  & 3.85  & 2.90  & 3.21  & 55.61 & 8.52 \\ \hline
    \end{tabularx}
\end{table}
We further compare the performance of the SAA and MIECP methods with varying $\theta$ and $\gamma$ under heterogeneous acceptance probabilities, as illustrated in \hyperref[fig:robust]{Figure \ref*{fig:robust}}. Specifically, \hyperref[fig:subrobusttheta]{Figure \ref*{fig:subrobusttheta}} reports the mean values as well as the 95\% confidence intervals (shaded areas) of models' objective values and the computing time ratios between the two models when $\theta$ varies, based on the D20–S80 instance. \hyperref[fig:subrobustgamma]{Figure \ref*{fig:subrobustgamma}} presents the results obtained by varying the number of supplies while fixing $\theta = 4$, the number of demands at 20. Overall, the MIECP method requires less than 3\% of the average computation time of the SAA method, and the difference in objective values between MIECP and SAA is within 0.5\%. Moreover, when $\theta$ is relatively small or $\gamma$ is large, MIECP can outperform SAA in terms of solution quality. 

\subsection{Case Study: Full Truck Alliance}\label{sec:casestudy}
We conduct a case study using a large-scale data set from Full Truck Alliance (\url{https://www.fulltruckalliance.com/en}), which contains platform-level information on more than 100,000 drivers and shippers, and over one million historical orders in 2018 \citep{competation}. We first describe the calibration procedure for the main parameters.
We set the maximum number of recommendations $\theta$ to 5 because the ratio of drivers to shippers from historical assignments falls within the range between 4 and 6. From the platform’s perspective, the utility function needs to account for the characteristics of the driver, the order, and the platform itself, as detailed in \hyperref[tab:utility]{Table \ref*{tab:utility}}. Note that due to the negative impact of distance on utility, $\overline{d}_{ij}$ is defined as the maximum distance minus the actual distance, so that shorter distances correspond to higher utility. For the driver's familiarity, we calculate it based on the correlation between the current order’s origin–destination pair and the driver’s frequently traveled routes. 

\begin{figure}[!htbp]
      \centering
      \includegraphics[width=0.6\linewidth]{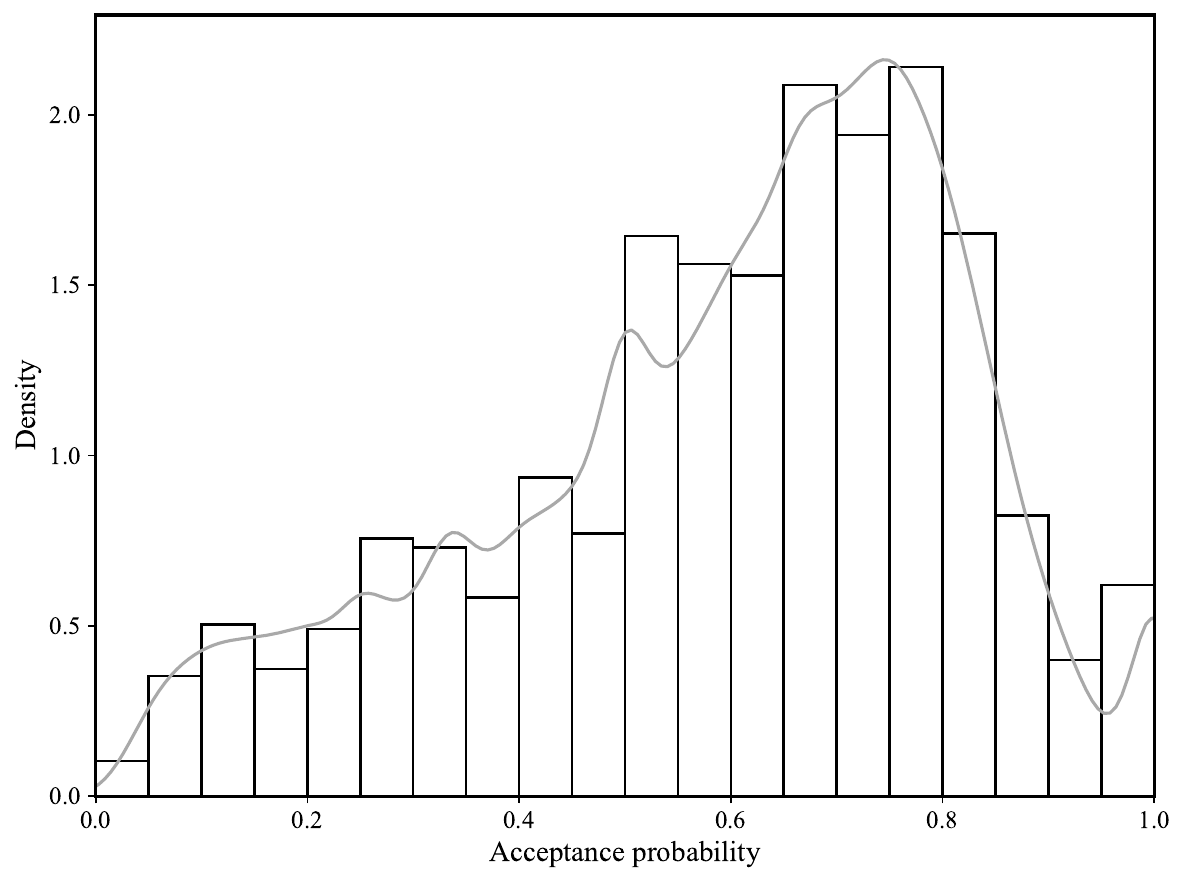}
      \caption{Distribution of the driver acceptance probability}
      \label{fig:distributionp}
    \end{figure}

Because the data did not directly reveal the driver's detailed acceptance response to every job, we consider a parsimonious linear utility function as follows. Given the greater importance of distance and familiarity in practical matching, we set their coefficients in the utility function to 0.3, the coefficient for order revenue to 0.2, and the coefficients for the driver’s and shipper’s evaluation to 0.1. We further normalize each attribute to the range of $[0,1]$, and the resulting utility function is assumed to take the form of $u_{ij}=0.1e^{D}_{j}+0.1e^{O}_{i}+0.2r^{O}_{i}+0.3\overline{d}_{ij}+0.3f_{ij}$. Robustness checks on the utility coefficients are provided at the end of this section and show that the main conclusion does not rely on these specifics. Moreover, the distribution of historical average order acceptance probabilities for over 100,000 drivers ($p_{j}^{H}$) is shown in \hyperref[fig:distributionp]{Figure \ref*{fig:distributionp}}, with the 25th, 50th, and 75th percentiles being 0.45, 0.64, and 0.76, respectively (acceptance probabilities of 0.05 and 0.95 correspond to the 0.5th and 96.9th percentiles, respectively). To better capture driver heterogeneity, we assume that the probability of the driver accepting the order is related to the driver's historical average acceptance rate $p_{j}^{H}$, the relative distance from the current location to the origin of the order $\overline{d}_{ij}^{R}$ (normalized to the range of $[-1,1]$), and the relative familiarity of the order route $f_{ij}^{R}$ (normalized to the range of $[-1,1]$). Taking into account these three factors, the acceptance probability function follows $p_{ij}=p_{j}^{H}+0.025\overline{d}_{ij}^{R}+0.025f_{ij}^{R}$, when $p_{j}^{H}\in [0.05,0.95]$, and otherwise, $p_{ij}=p_{j}^{H}$. 
\begin{table}[htbp]
    \caption{Attribute definition in the platform utility function}
    \label{tab:utility}
    \centering
    \renewcommand{\arraystretch}{1.15} 
    \begin{tabularx}{\textwidth}{c c >{\raggedright\arraybackslash}X}
    \toprule
    Class & Attribute & Definition \\ \midrule
    Driver & $e^{D}_{j}$ & Driver evaluation score \\ \midrule
    \multirow{2}{*}{Order} 
     & $e^{O}_{i}$ & Shipper (customer) evaluation score \\
     & $r^{O}_{i}$ & Order revenue \\ \midrule
    \multirow{2}{*}{Platform} 
     & $\overline{d}_{ij}$ & Maximum distance minus the distance between the driver's current location and the origin of the order \\
     & $f_{ij}$ & Driver's familiarity with the route of order \\ \bottomrule
    \end{tabularx}
\end{table}

Based on the above setup, we generate instances by sampling orders and drivers from the historical samples. Note that in long-haul freight transportation, shippers (demands) often impose hard constraints on vehicle types, while each driver (supply) typically has only one vehicle type. Thus, we sample three instances from drivers and orders of a single vehicle type, with sizes D100-S450 ($\gamma < \theta$), D100-S500 ($\gamma = \theta$), and D100-S550 ($\gamma > \theta$), to simulate different values of $\gamma$. 
In addition to DAP, DIMP, SAA, and MIECP, we also evaluate a commonly used industry benchmark called the nearby-priority policy (NPP) that recommends the order to the drivers who are close to its origin, and assigns each order to the nearest driver who accepts the recommendation \citep{industry} (the detailed formulation is presented in \ref{sec:anpp}). We compare their performance over 100 samples of driver responses, where each driver’s response is randomly generated based on the given acceptance probability. The evaluation results are summarized in \hyperref[tab:case]{Table \ref*{tab:case}}, where the reported values represent the relative gap (\%) with the best-performing solution. Overall, MIECP consistently obtains the best matching solution among the tested methods. The improvement is particularly notable when $\gamma < \theta$. 

{
\begin{table}[htbp]
        \caption{Performance evaluation (gap with the best solution found) using Full Truck Alliance data}
        \label{tab:case}
        \centering
\begin{tabular}{cccc}
\hline
 Model     & $\gamma < \theta$ & $\gamma = \theta$ & $\gamma > \theta$ \\ \hline
DAP   & 6.12              & 0.21              & 0.15              \\
DIMP  & 9.60              & 9.68              & 8.66              \\
NPP   & 4.80              & 2.00              & 2.09              \\
SAA   & 0.05              & 0.28              & 0.09              \\
MIECP & 0.00              & 0.00              & 0.00              \\ \hline
\end{tabular}
    \end{table}
}

\paragraph{Performance with misspecified acceptance probabilities.} Moreover, to understand the performance of matching policies under acceptance probability misspecifications, we consider four types of out-of-sample scenarios: randomly lower probabilities (`Out-L'); randomly higher probabilities (`Out-H'); small random perturbations around in-sample probabilities (`Out-NS'); large random perturbations around in-sample probabilities (`Out-NL'). Similarly to the in-sample tests, drivers’ acceptance behaviors follow Bernoulli distributions according to the specified acceptance probabilities as follows:
\begin{itemize}
    \item \textit{Scenario `Out-L'}: Randomly generated within $[\max(0,p_{ij}-0.05),p_{ij}]$
    \item \textit{Scenario `Out-H'}: Randomly generated within $[p_{ij},\min(p_{ij}+0.05,1)]$
    \item \textit{Scenario `Out-NS'}: Randomly generated within $[\max(0,p_{ij}-0.025),\min(p_{ij}+0.025,1)]$
    \item \textit{Scenario `Out-NL'}: Randomly generated within $[\max(0,p_{ij}-0.1),\min(p_{ij}+0.1,1)]$.
\end{itemize}

The results are summarized in \hyperref[tab:caseout]{Table \ref*{tab:caseout}}. We observe that MIECP achieves the best out-of-sample performance in most cases, except for `Out-H'. Consistent with the in-sample results, the improvement of the MIECP model is most significant when $\gamma < \theta$. These results suggest that the MIECP approximation is not only preferable in computational efficiency but also robust against acceptance probability misspecification (while the performance of SAA is more volatile).
\begin{table}[htbp]
\caption{Out-of-sample performance evaluation using Full Truck Alliance data}
\label{tab:caseout}
\centering
\begin{adjustbox}{width=\textwidth}
\begin{tabular}{c rrrr c rrrr c rrrr}
\hline
\multirow{2}{*}{\rule{0pt}{2.6ex}Method} 
& \multicolumn{4}{c}{$\gamma < \theta$} & 
& \multicolumn{4}{c}{$\gamma = \theta$} & 
& \multicolumn{4}{c}{$\gamma > \theta$} \\[0.4ex]
\cline{2-5} \cline{7-10} \cline{12-15}
& \rule{0pt}{2.4ex}Out-L & Out-H & Out-NS & Out-NL 
& & Out-L & Out-H & Out-NS & Out-NL 
& & Out-L & Out-H & Out-NS & Out-NL \\ 
\hline
DAP   & 6.07  & 6.06 & 6.00  & 5.85  & & 0.31  & 0.36 & 0.15 & 0.21  & & 0.10  & 0.28 & 0.00 & 0.19 \\
DIMP  & 11.69 & 7.56 & 11.07 & 32.07 & & 12.16 & 6.53 & 9.71 & 31.07 & & 11.37 & 6.52 & 9.02 & 27.90 \\
NPP   & 4.82  & 3.57 & 4.56  & 4.58  & & 1.98  & 2.00 & 2.02 & 2.17  & & 2.10  & 2.00 & 2.06 & 1.96 \\
SAA   & 0.36  & 0.00 & 0.13  & 3.48  & & 0.71  & 0.00 & 0.35 & 4.54  & & 0.77  & 0.00 & 0.35 & 3.32 \\
MIECP & 0.00  & 0.18 & 0.00  & 0.00  & & 0.00  & 0.25 & 0.00 & 0.00  & & 0.00  & 0.24 & 0.00 & 0.00 \\
\hline
\end{tabular}
\end{adjustbox}
\end{table}



\paragraph{Additional robustness checks.} To demonstrate that variations in the coefficient of the utility function and the acceptance probability function have no significant impact on the performance of our model, we conduct robustness checks based on three additional types of instances (the instance size is the same as the case study). Specifically, we randomly generate ten utility instances for each, where the coefficients of each term in the utility function are uniformly drawn from $[0.1, 0.3]$. Furthermore, we consider three types of acceptance probability function coefficients ($f(p)$): (i) $f(p)^{1}$: consistent with the aforementioned description. (ii) $f(p)^{2}$: $p_{ij}=p_{j}^{H}+0.05\overline{d}_{ij}^{R}+0.05f_{ij}^{R}$, when $p_{j}^{H}\in [0.1,0.9]$, and otherwise, $p_{ij}=p_{j}^{H}$. (iii) $f(p)^{3}$: $p_{ij}=p_{j}^{H}$. \hyperref[tab:caserobust]{Table \ref*{tab:caserobust}} reports the relative gap (\%) of the four methods with the best solution found, across ten randomly generated utility function coefficient settings for each instance. We observe that MIECP consistently delivers the best performance, with a relative average gap of less than 0.01\%, whereas DAP, DIMP, NPP, and SAA have average gaps of 3.06\%, 16.42\%, 7.73\%, and 0.72\%, respectively. Note that, compared with the main case study, the average performance of NPP deteriorates, mainly due to the uncertainty in the utility function coefficients (particularly when the weight of distance decreases). 

\begin{table}[htbp]
    \caption{Performance evaluation with varying utility functions and the acceptance probability functions}
    \label{tab:caserobust}
    \centering
    \begin{tabular}{c rrr c rrr c rrr}
    \hline
    \multirow{2}{*}{\rule{0pt}{2.6ex}Method} 
    & \multicolumn{3}{c}{$\gamma < \theta$} & 
    & \multicolumn{3}{c}{$\gamma = \theta$} & 
    & \multicolumn{3}{c}{$\gamma > \theta$} \\[0.4ex]
    \cline{2-4} \cline{6-8} \cline{10-12}
    & \rule{0pt}{2.4ex}$f(p)^1$ & $f(p)^2$ & $f(p)^3$ 
    & & $f(p)^1$ & $f(p)^2$ & $f(p)^3$ 
    & & $f(p)^1$ & $f(p)^2$ & $f(p)^3$ \\ 
    \hline
    DAP   & 7.00  & 6.92  & 7.13  & & 1.35  & 1.52  & 1.09  & & 0.82  & 0.87  & 0.83 \\
    DIMP  & 16.10 & 14.34 & 17.50 & & 17.30 & 15.70 & 18.56 & & 16.20 & 14.47 & 17.66 \\
    NPP   & 8.41  & 8.80  & 8.08  & & 7.83  & 8.02  & 7.68  & & 6.92  & 7.07  & 6.82 \\
    SAA   & 0.44  & 0.39  & 0.35  & & 0.91  & 0.90  & 0.93  & & 0.86  & 0.73  & 1.02 \\
    MIECP & 0.04  & 0.02  & 0.00  & & 0.00  & 0.00  & 0.00  & & 0.00  & 0.00  & 0.00 \\
    \hline
    \end{tabular}
\end{table}


\section{Concluding Remarks}
We study a two-stage matching strategy that is vital for freight-matching and other crowdsourcing platforms to hedge against uncertain supply acceptance responses when facing a large supply pool. The acceptance uncertainty introduces nonlinearity that makes conventional stochastic programming approaches computationally expensive, particularly for large-scale applications. Although deterministic approximation approaches enjoy scalable LP formulations, we construct cases to show that their performance can be arbitrarily bad. This motivates us to develop a principled convex approximation approach, which yields an MIECP formulation that can be efficiently solved. We prove worst-case performance guarantees and establish conditions where the proposed approximation strictly dominates a deterministic approximation approach. Numerical results support the superior performance of the ECP-based approximation in terms of both solution quality and time.

While we focus on a one-shot matching problem in this study (for cases where the decision frequency is low and/or the jobs are not very time sensitive), the proposed approximation policy based on MIECP could be used to develop advanced dynamic matching algorithms (e.g., using value function approximations), and our performance analysis could be an important building block to establishing performance guarantees of dynamic policies. For example, in a food donation recommendation problem studied by \cite{lee2025offer}, the donation jobs are recommended/offered to agencies sequentially across multiple stages and may be rejected. Our approach offers a promising means to optimizing the corresponding offer scheduling problem, where assignment of donations to accepted agencies is endogenized. 

Because the considered matching problem has a generic structure shared by other operational problems, and the inner maximization with uncertain acceptance is reminiscent of recourse decisions in the face of arc disruptions (of a network), the proposed MIECP formulation could be leveraged to solve supply chain resilience planning problems. In particular, one may extend our framework and analysis to general graphs with more than two types of agents, where the arcs are susceptible to disruptions. Moreover, when matching decisions are high-stakes, the failure of matching should be more carefully hedged against, and one may develop robust optimization approaches to design risk-averse matching policies.

\OneAndAHalfSpacedXI
\bibliographystyle{informs2014trsc}
\bibliography{mybib}

\ECSwitch

\counterwithin{table}{section}
\counterwithin{figure}{section}
\numberwithin{equation}{section}  
\makeatletter

\renewcommand{\p@subfigure}{} 
\renewcommand{\thesubfigure}{\thesection.\arabic{figure}(\alph{subfigure})}
\makeatother
\vspace{-20pt}

\section{Technical Proofs} \label{sec:proof}
\subsection{Proof of Lemma \ref{pro1}}
\textit{Proof.} Given any feasible recommendation solution of SP, it is also feasible to R-SP, and vice versa. Then, if the objective value of SP and R-SP are the same for any recommendation solution, then the optimal objective value and the set of optimal solutions of model SP will coincide with those of R-SP. Consider a fixed recommendation set for demand $i$, denoted by $\{j_1, j_2, \dots, j_k\}$ such that $x_{ij}=1$ for these $j$, and sorted in descending order of utility, i.e., $u_{ij_1} \ge u_{ij_2} \ge \dots \ge u_{ij_k}$, the expected utility contributed by demand $i$ is
\begin{equation} \label{eq:prob}
    \mathbb{E}_\xi\left[ \max_{j} u_{ij} \cdot \tilde{\xi}_{ij} \right] = \sum_{r=1}^k p(1-p)^{r-1} u_{ij_r} = \sum_{r=1}^k p_r u_{ij_r},
\end{equation}
where $p_r$, by definition, equals the probability that the $r$-th ranked supply is the first one to accept. We then prove that given a fixed recommendation set, the optimal solution $w_{ir}$ from R-SP is equal to $u_{ij_r}$ by contradiction. Suppose that according to the optimal solution $\mathbf{x}_R$ of R-SP, there exists $r_1<r_2$ such that $w_{ir_1} < w_{ir_2}$, that is, $w_{ir}$'s for demand $i$ are not sorted in descending order of utility. Then we can construct a feasible solution  $\mathbf{x}^{\prime}_R$ by exchanging $r_1$ and $r_2$ of demand $i$, while keeping all other variables unchanged. Because $p_{r_1}>p_{r_2}$ (by definition of $p_r$), we have $p_{r_1}w_{ir_1}+p_{r_2}w_{ir_2}<p_{r_1}w_{ir_2}+p_{r_2}w_{ir_1}$, implying that the objective value of $\mathbf{x}^{\prime}_R$ is greater than that of $\mathbf{x}_R$, which contradicts our assumption. 
Therefore, under $\mathbf{x}_R$, $w_{ir}$'s must be sorted in descending order and take values in $\{u_{ij_r}\}$ due to constraints \eqref{eq:r8}-\eqref{eq:r9} for all $i\in\mathcal{A}^D$. It follows that SP and R-SP yield the same expected utility, thereby proving the proposition.\halmos
\subsection{Proof of Corollary \ref{corollary0}}
{\color{black}
\textit{Proof.} We first introduce binary variables $y_{ijr}$ indicating whether supply $j$ is the $r$th largest-utility element in demand $i$'s recommendation set. In this way, the objective function becomes: $$ \max \  \sum_{i \in \mathcal{A}^{D}}\sum_{j \in \mathcal{A}^S}\sum_{r \in \mathcal{R}}p_{r}\cdot u_{ij}\cdot y_{ijr}.$$ Since the binary variable $z_{ijr}$ indicates whether supply $j$ is among the top $r$ utility supplies in demand $i$'s recommendations, variables $z_{ijr}$ and $y_{ijr}$ satisfy the following relationship (since $z_{ijr}\ge z_{ij(r-1)}$, there are a
 total of five possible cases):
\begin{align*}
   & \text{if } z_{ijr}=1 \text{ and } z_{ij(r-1)}=0 \text{, then } y_{ijr}=1, \forall i \in \mathcal{A}^{D}, j \in \mathcal{A}^S, r \in \mathcal{R}\setminus\{1\}, \\
   &  \text{if } z_{ijr}=1 \text{ and } z_{ij(r-1)}=1 \text{, then } y_{ijr}=0, \forall i \in \mathcal{A}^{D}, j \in \mathcal{A}^S, r \in \mathcal{R}\setminus\{1\}, \\
   &  \text{if } z_{ijr}=0 \text{ and } z_{ij(r-1)}=0 \text{, then } y_{ijr}=0, \forall i \in \mathcal{A}^{D}, j \in \mathcal{A}^S, r \in \mathcal{R}\setminus\{1\}, \\
    & \text{if } z_{ij1}=1  \text{, then } y_{ij1}=1, \forall i \in \mathcal{A}^{D}, j \in \mathcal{A}^S, \\
   & \text{if } z_{ij1}=0  \text{, then } y_{ij1}=0, \forall i \in \mathcal{A}^{D}, j \in \mathcal{A}^S.
\end{align*}
These logical constraints can be linearized as follows:
\begin{subequations}
\begin{align}
    &  z_{ijr}-z_{ij(r-1)}=y_{ijr}, \forall i \in \mathcal{A}^{D}, j \in \mathcal{A}^S, r \in \mathcal{R}\setminus\{1\}, \label{eq:linear1}\\
    & z_{ij1}=y_{ij1},\forall i \in \mathcal{A}^{D}, j \in \mathcal{A}^S.\label{eq:linear2}
\end{align}
\end{subequations}
For a given $(i,j)$ pair, by summing the above constraints, we obtain:
\begin{align*}
   z_{ij1}+ \sum_{r \in \mathcal{R}\setminus\{1\}}\left(z_{ijr}-z_{ij(r-1)}\right)=z_{ij|\mathcal{R}|} = \sum_{r \in \mathcal{R}}y_{ijr}.
\end{align*}
Moreover, based on the definition of $\mathbf{y}$ variables, for each demand-rank pair, at most one supply can be assigned, thus we impose the following constraints:
$$\sum_{j\in \mathcal{A}^{S}}y_{ijr} \le 1, \forall i \in \mathcal{A}^{D}, r \in \mathcal{R}.$$
Based on the definition of $z_{ijr}$, we have $z_{ij|\mathcal{R}|}=x_{ij}$, then we have $\sum_{r \in \mathcal{R}}y_{ijr}=x_{ij}$. Note that, By replacing $x_{ij}$ with $\sum_{r \in \mathcal{R}}y_{ijr}$ and replacing $w_{ir}$ with $u_{ij}y_{ijr}$, the R-SP model can be equivalently reformulated as:
\begin{subequations}
    \begin{align}
        \textbf{[R-SP2]}: \max \  &\sum_{i \in \mathcal{A}^{D}}\sum_{j \in \mathcal{A}^S}\sum_{r \in \mathcal{R}}p_{r}\cdot u_{ij}\cdot y_{ijr}\\
         {\rm s.t.} \ & \hyperref[eq:r4]{(\ref*{eq:r4})},\hyperref[eq:r5]{(\ref*{eq:r5})},\hyperref[eq:r7]{(\ref*{eq:r7})},\hyperref[eq:linear1]{(\ref*{eq:linear1})}, \hyperref[eq:linear2]{(\ref*{eq:linear2})}, \nonumber  \\
          &\sum_{i \in \mathcal{A}^{D}}\sum_{r \in \mathcal{R}}y_{ijr} \le 1, \forall j\in \mathcal{A}^{S},\label{eq:lp1}\\
          &\sum_{j\in \mathcal{A}^{S}}y_{ijr} \le 1, \forall i \in \mathcal{A}^{D}, r \in \mathcal{R}, \label{eq:lp2}\\
        & \sum_{j\in \mathcal{A}^{S}}\sum_{r \in \mathcal{R}}y_{ijr}  \le \theta , \forall i \in \mathcal{A}^{D}, \label{eq:lp3}\\
         & z_{ijr} \le \sum_{r^{\prime} \in \mathcal{R}}y_{ijr^{\prime}}, \forall i \in \mathcal{A}^{D}, j\in \mathcal{A}^{S}, r \in \mathcal{R}, \\
          &  t_{ir} \ge \sum_{j\in \mathcal{A}^{S}}\sum_{r^{\prime} \in \mathcal{R}}y_{ijr^{\prime}}-\theta v_{ir}, \forall i \in \mathcal{A}^{D}, r \in \mathcal{R}\\
         & u_{ij}\cdot y_{ijr} \le u_{ij}+ M^{1}_{i}(1-z_{ijr}),\forall i \in \mathcal{A}^{D}, j\in \mathcal{A}^{S}, r \in \mathcal{R},\\
         & u_{ij}\cdot y_{ijr}\le M^{1}_{i}v_{ir}, \forall i \in \mathcal{A}^{D}, r \in \mathcal{R},\\
         & y_{ijr},z_{ijr},v_{ir} \in \{0,1\}, t_{ir} \in \mathbb{N},  \forall i \in \mathcal{A}^{D}, j\in \mathcal{A}^{S},  r \in \mathcal{R}.
    \end{align}
\end{subequations}

LP formulation of \citet{ekbatani2026lyft1} has the same objective as the R-SP2 model, but retains only constraints \eqref{eq:lp1} and \eqref{eq:lp2} (binary variables $\mathbf{y}$ can be relaxed to continuous since the coefficient matrix is totally unimodular), which is denoted as model R-SP-LP. A minor difference is that they only require the size of the set $\mathcal{R}$ to be no greater than the number of supplies, whereas we impose an explicit upper bound of $\theta$. 
We next prove that {R-SP-LP} is an exact projection of {R-SP2} onto the $\mathbf{y}$-space; consequently, the two models attain the same optimal objective value (and admit the same set of optimal $\mathbf{y}$ solutions). Given any feasible solution $\mathbf{\bar{y}}$ of R-SP-LP, constraints \eqref{eq:lp3} are always satisfied, since it can be obtained by summing constraints \eqref{eq:lp2} over all $r\in \mathcal{R}$. For the remaining variables, there always exists a feasible solution in R-SP2 while keeping $\mathbf{y}$ fixed. Intuitively, we can construct a feasible solution: $\bar{z}_{ijr}=\sum_{k=1}^{r}\bar{y}_{ijk}$, $\bar{t}_{ir}=\min(\sum_{j\in \mathcal{A}^{S}}\sum_{r^{\prime}\in \mathcal{R}}\bar{y}_{ijr^{\prime}},r)$ and $\bar{v}_{ir}= \mathbbm{1}(\sum_{j\in \mathcal{A}^{S}}\sum_{r^{\prime}\in \mathcal{R}}\bar{y}_{ijr^{\prime}}\ge r)$.

Therefore, $\bar{\mathbf y}$ is extendable to a feasible solution of {R-SP2}. This implies that the projection of the feasible set of {R-SP2} onto the $\mathbf y$-space coincides with the feasible set of {R-SP-LP}. Since the two models share the same objective function in $\mathbf y$, they attain the same optimal objective value and share the same set of optimal $\mathbf y$ solutions. \halmos}
\subsection{Proof of Proposition \ref{prop:unbound}}
\textit{Proof.} 
Suppose there exists a set of $\lceil \gamma |\mathcal{A}^{D}|/\theta \rceil$ demands, denoted as $\mathcal{A}^{D\prime}$, such that $u_{i^{\star}j} = b$ for $i^{\star}\in \mathcal{A}^{D\prime}$ and $j\in \mathcal{A}^S$. For all remaining demand-supply pairs, the utility is equal to $a$. The acceptance probabilities are all equal to $p$. In this case, under the DAP policy, all supplies will be recommended demands in the set $\mathcal{A}^{D\prime}$, and the resulting expected objective value $Obj^{E} (\mathbf{x}^{D})$ will be at most $\lceil \gamma |\mathcal{A}^{D}|/\theta \rceil b$ (each supply can be matched at most one demand). Meanwhile, the value of $Obj^{E} (\mathbf{x}^{\star})$ is no less than $\min \left(|\mathcal{A}^{D}|,\gamma|\mathcal{A}^{D}|\right)p\cdot a$. Therefore, in this case, we have:
\begin{align*}
    \frac{Obj^{E} (\mathbf{x}^{\star})-Obj^{E} (\mathbf{x}^{D})}{Obj^{E} (\mathbf{x}^{\star})} > 1-\frac{\lceil \gamma |\mathcal{A}^{D}|/\theta \rceil b}{\min \left(|\mathcal{A}^{D}|,\gamma|\mathcal{A}^{D}|\right)p\cdot a}.
\end{align*}
Specifically, when the proportion $\gamma / \theta$ approaches zero, the approximation gap will approach 1. \Halmos
\subsection{Proof of Proposition \ref{dapnew}}
\textit{Proof.} Under the condition of $\gamma \ge \theta >1$, it is intuitive to find the value of $Obj^{E} (\mathbf{x}^{M})$ is not greater than $|\mathcal{A}^D|\cdot p\cdot b$ and the value of $Obj^{E} (\mathbf{x}^{\star})$ is no less than $|\mathcal{A}^D|\cdot \left(1-\left(1-p\right)^\theta \right)\cdot a$.
Therefore, in this case, we have:
\begin{align*}
    \frac{Obj^{E} (\mathbf{x}^{\star})-Obj^{E} (\mathbf{x}^{M})}{Obj^{E} (\mathbf{x}^{\star})} > 1-\frac{p\cdot b}{\left(1-\left(1-p\right)^\theta \right)\cdot a}.
\end{align*}
Specifically, when $b\rightarrow a$, $p \rightarrow 0$ and $\theta \rightarrow \infty$, the approximation gap will approach 1. \Halmos
\subsection{Proof of Corollary \ref{corollary1}} 
\textit{Proof.} Given solution $x$, we have for $i\in \mathcal{A}^D$
\begin{align*}
     \mathbb{E}_{\mathbb{P}} \left [\max_{j\in \mathcal{A}^{S}} u_{ij} \tilde{\xi}_{ij} x_{ij}\right] &\leq \tau \mathbb{E}_{\mathbb{P}} \left [\log \left(\sum_{j\in \mathcal{A}^{S},x_{ij}=1}\exp{\left(\frac{u_{ij} \tilde{\xi}_{ij}}{\tau}\right)}\right)\right] & \quad (\text{By \hyperref[lemma0]{Lemma \ref*{lemma0}}}) \\
    & \leq \tau \log \left[\mathbb{E}_{\mathbb{P}}\left( \sum_{j\in \mathcal{A}^{S}, x_{ij}=1}\exp{\left(\frac{u_{ij} \tilde{\xi}_{ij}}{\tau}\right)}\right)\right] & \quad (\text{By Jensen's inequality}) \\
    & = \tau \log \left[\sum_{j\in \mathcal{A}^{S}, x_{ij}=1}\mathbb{E}_{\mathbb{P}} \left(\exp{\left(\frac{u_{ij} \tilde{\xi}_{ij}}{\tau}\right)}\right)\right] &\quad (\text{Linearity of expectation})\\
    & = \tau \log \left[\sum_{j\in \mathcal{A}^{S}, x_{ij}=1}\left(\exp{\left(\frac{u_{ij} }{\tau}\right)}p_{ij} + 1\cdot (1-p_{ij})\right)\right] & \quad (\text{Bernoulli distribution of } \tilde{\xi}) \\
    &= \tau \log \left[\sum_{j\in \mathcal{A}^{S}}\left((\exp{\left(\frac{u_{ij} }{\tau}\right)}-1)p_{ij} + 1\right)x_{ij}\right].  & \quad (\text{Substituting into variable } \mathbf{x})
\end{align*}
Summing the above inequalities over all demands leads to the desired result. \Halmos
\subsection{Proof of Proposition \ref{prop:not_infty_minimizer}}
{\color{black}
\textit{Proof.} We prove in four steps. First, we prove the convexity of $f(\tau)$. In steps 2 and 3, we derive the limits of $f_i'(\tau)$ as $\tau \to \infty$ and $\tau \to 0^+$. Finally, we discuss the sufficient conditions for $f_i'(0^+)<0$ and $f_i'(\infty)>0$.

\noindent\textit{Step 1: Convexity of $f(\tau)$.}
Fix a feasible recommendation solution $\mathbf{x}$ and we define $f(\tau)=\sum_i f_i(\tau)$ with
$$f_i(\tau)=\tau\log S_i(\tau),\qquad S_i(\tau)=\sum_{j\in A_i}\Bigl((1-p_{ij})+p_{ij}\exp({u_{ij}/\tau})\Bigr).$$
We can obtain the second derivative of $f_i(\tau)$ as
\[
f_i''(\tau)=
\frac{1}{\tau^3}\left(
\frac{\sum_{j\in A_i}p_{ij}u_{ij}^2e^{u_{ij}/\tau}}{S_i(\tau)}
-\left(\frac{\sum_{j\in A_i}p_{ij}u_{ij}e^{u_{ij}/\tau}}{S_i(\tau)}\right)^2
\right),\qquad \tau>0.
\]
Introducing weights $w_{ij}(\tau):=\frac{p_{ij}e^{u_{ij}/\tau}}{S_i(\tau)}\ge 0$, then we have
\[
\sum_{j\in A_i} w_{ij}(\tau)u_{ij}^2-\Bigl(\sum_{j\in A_i} w_{ij}(\tau)u_{ij}\Bigr)^2
=\frac12\sum_{j\in A_i}\sum_{k\in A_i} w_{ij}(\tau)w_{ik}(\tau)\,(u_{ij}-u_{ik})^2\ge 0,
\]
which implies $f_i''(\tau)\ge 0$ for all $\tau>0$. Therefore, $f_i(\tau)$ is convex on $(0,\infty)$, for each $i\in \mathcal{A}^D$. Since $f_i(\tau)=\tau\log S_i(\tau)$ is convex on $(0,\infty)$, it follows that $f(\tau)$ is convex and $f'(\tau)$ is nondecreasing. 
Let
\[
N_i(\tau):=\sum_{j\in A_i}p_{ij}u_{ij}\exp({u_{ij}/\tau}).
\]
A direct differentiation yields, for $\tau>0$,
\begin{equation*}
f_i'(\tau)=\log S_i(\tau)-\frac{1}{\tau}\frac{N_i(\tau)}{S_i(\tau)},\qquad 
f'(\tau)=\sum_i f_i'(\tau).
\end{equation*}
Since we cannot obtain a closed-form solution for the stationary points, we discuss the function values at the endpoints.

\noindent\textit{Step 2: Endpoint derivative at $\tau\to\infty$.} Using $\exp(u_{ij}/\tau)=1+\frac{u_{ij}}{\tau}+O(\tau^{-2})$, we obtain
\[
S_i(\tau)=|A_i|+\frac{1}{\tau}\sum_{j\in A_i}p_{ij}u_{ij}+O(\tau^{-2}), \quad N_i(\tau)=\sum_{j\in A_i}p_{ij}u_{ij}+\frac{1}{\tau}\sum_{j\in A_i}p_{ij}u_{ij}^2+O(\tau^{-2}).
\]

Consequently, we have
\[
\frac{N_i(\tau)}{S_i(\tau)}=\frac{1}{|A_i|}\sum_{j\in A_i}p_{ij}u_{ij}+O(\tau^{-1}),
\qquad
\frac{1}{\tau}\frac{N_i(\tau)}{S_i(\tau)}\to 0,
\]
and thus
\[
\lim_{\tau\to\infty} f_i'(\tau)=\log|A_i|,
\qquad
\lim_{\tau\to\infty} f'(\tau)=\sum_i \log|A_i|.
\]

\noindent\textit{Step 3: Endpoint derivative at $\tau\to 0^+$.} We denote by $u_i^{\max}$ the maximum utility among the suppliers recommended to demand $i$, and let $J_i^\star$ be the corresponding set of suppliers attaining this maximum (under \hyperref[assu1]{Assumption \ref*{assu1}}, this set contains only one element.) For any recommended $j\notin J_i^\star$, since $u_{ij}<u_i^{\max}$, we have
$\exp(u_{ij}/\tau)/\exp(u_i^{\max}/\tau)\to 0$.
Hence
\[
S_i(\tau)=\sum_{j\in A_i}(1-p_{ij})+\exp(u_i^{\max}/\tau)\Bigl(\sum_{j\in J_i^\star}p_{ij}+o(1)\Bigr),
\]
and the exponential term dominates the constant term as $\tau\to 0^+$, yielding
\[
S_i(\tau)=\exp(u_i^{\max}/\tau)\Bigl(\sum_{j\in J_i^\star}p_{ij}+o(1)\Bigr).
\]
Therefore,
\[
\log S_i(\tau)=\frac{u_i^{\max}}{\tau}+\log\Bigl(\sum_{j\in J_i^\star}p_{ij}\Bigr)+o(1).
\]
Similarly,
\[
N_i(\tau)=\exp(u_i^{\max}/\tau)\Bigl(u_i^{\max}\sum_{j\in J_i^\star}p_{ij}+o(1)\Bigr),
\qquad
\frac{N_i(\tau)}{S_i(\tau)}=u_i^{\max}+o(1).
\]
Substituting into $f_i'(\tau)$ gives
\[
\lim_{\tau\to 0^+} f_i'(\tau)=\log\Bigl(\sum_{j\in J_i^\star}p_{ij}\Bigr),
\qquad
\lim_{\tau\to 0^+} f'(\tau)=\sum_i \log\Bigl(\sum_{j\in J_i^\star}p_{ij}\Bigr).\]

\noindent\textit{Step 4: Sign of the endpoint derivative.}
If $|A_i|\ge 1 $ for all demands and there exists at least one demand $i_0$ such that $|A_{i_0}|\ge 2$, the value of $\lim_{\tau\to\infty} f'(\tau)=\sum_i \log|A_i|$ is positive. Moreover, under \hyperref[assu1]{Assumption \ref*{assu1}}, $J_i^\star$ contains one element and $p_{ij}<1$, we have $\lim_{\tau\to 0^+} f'(\tau)=\sum_i \log\Bigl(\sum_{j\in J_i^\star}p_{ij}\Bigr)<0$. Then, for convex function $f(\tau)$, the minimizer ($f'(\tau^{\star})=0$) exists and lies in $(0,\infty)$.\Halmos}

\subsection{Proof of Proposition \ref{pro:miecp}}
\textit{Proof.} We first denote the inner linear term of objective \hyperref[eq:5a]{(\ref*{eq:5a})} $\sum_{j \in \mathcal{A}^{S} }\exp{\left(\frac{u_{ij}}{\tau}\right)}p_{ij}x_{ij}$ as $Y_{i}$. Then, we introduce an auxiliary variable $Z_{i}$ for each $i \in \mathcal{A}^{D}$, representing the logarithmic term:
\[
    Z_{i}=\log(Y_{i}), \forall i \in \mathcal{A}^{D}.
\]
Because the objective is to maximize the summation of $Z_i$, the above equality constraint can be written as an inequality constraint $Z_{i} \le \log(Y_{i})$, which can then be modeled using an exponential cone $\left( Y_{i},1, Z_{i}\right) \in \mathcal{K}_{\exp}$, resulting in a mixed-integer exponential cone program (MIECP). \Halmos
\subsection{Proof of Lemma \ref{lemma:bounds}}
\textit{Proof.} We prove each of the four inequalities as follows. 

\noindent (i) \textit{Lower bound of $V_{i}^{E}(\mathbf{x})$.} Let $j^{\star}$ be the index such that $u_{ij^{\star}} = U_i(x)$ (for $x_{ij^{\star}}=1$). The expectation $V_i^{E}(\mathbf{x})$ is at least the expected value of any single term in the maximum:
$$ V_i^{E}(\mathbf{x}) = \mathbb{E}\left[\max_{j\in \mathcal{A}^{S},x_{ij}=1} u_{ij}\tilde{\xi}_{ij}\right] \ge \mathbb{E}[u_{ij^{\star}} \tilde{\xi}_{ij^{\star}}] = u_{ij^{\star}} p_{ij^{\star}} = U_i(\mathbf{x}) p_{ij^{\star}} \ge \underline{p}U_i(\mathbf{x}). $$
(ii) \textit{Upper bound of $V_{i}^{E}(\mathbf{x})$.} By the definition of $\overline{q}=\left(1-(1-\overline{p})^{\theta}\right)$, we have
$$ V_i^{E}(\mathbf{x}) = \mathbb{E}\left[\max_{j\in \mathcal{A}^{S},x_{ij}=1} u_{ij}\tilde{\xi}_{ij}\right] \le \mathbb{E}\left[ \max_{j\in \mathcal{A}^{S},x_{ij}=1} U_{i}(\mathbf{x})\tilde{\xi}_{ij}\right] = \left( 1- \Pi_{j\in \mathcal{A}^{S},x_{ij}=1}(1-p_{ij}) \right)U_{i}(\mathbf{x})\le  \overline{q}U_{i}(\mathbf{x}). $$
(iii) \textit{Lower bound of $V_{i}^{A}(\mathbf{x})$.} By the property \hyperref[eq:property]{(\ref*{eq:property})} of log-sum-exp function, we have the following lower bound of function $V_{i}^{A}(\mathbf{x})$:
\begin{align*}
   V_{i}^{A}(\mathbf{x}) &= \tau \log \left( \sum_{j \in \mathcal{A}^{S},x_{ij}=1} p_{ij} \exp{\left(\frac{u_{ij}}{\tau}\right)} \right) \\
   & \ge \tau \log(\underline{p}) + \tau \log \left( \sum_{j \in \mathcal{A}^{S},x_{ij}=1} \exp{\left(\frac{u_{ij}}{\tau}\right)} \right) &  (\text{Since } p_{ij} \ge \underline{p})\\ 
    &\ge \tau \log(\underline{p}) + U_i(\mathbf{x}). & (\text{By \hyperref[lemma0]{Lemma \ref*{lemma0}}})
\end{align*}
(iv) \textit{Upper bound of $V_{i}^{A}(\mathbf{x})$.} Similarly, the upper bound of $V_{i}^{A}(\mathbf{x})$ can be derived as
\begin{align*}
     V_{i}^{A}(\mathbf{x}) &= \tau \log \left( \sum_{j \in \mathcal{A}^{S},x_{ij}=1} p_{ij} \exp{\left(\frac{u_{ij}}{\tau}\right)} \right)  \\
     & \le \tau \log(\overline{p})+\tau \log \left(\sum_{j \in \mathcal{A}^{S},x_{ij}=1 }\exp{\left(\frac{u_{ij}}{\tau}\right)}\right) & (\text{Since } p_{ij} \le \overline{p}) \\ 
    &\le \tau \log(\overline{p})+ U_i(\mathbf{x}) + \tau \log (\theta).   &  (\text{By \hyperref[lemma0]{Lemma \ref*{lemma0}}}) 
\end{align*}
\Halmos
\subsection{Proof of Theorem \ref{theoemhomo}}
\textit{Proof.} We prove the bound in four steps. In steps 1 and 2, we establish an upper bound of $Obj^{E} (\mathbf{x}^{\star})-Obj^{E}(\mathbf{x}^{A})$ for any realization of $u_{ij}$. Then, we bound $\mathbb{E}\left[{Obj^{E} (\mathbf{x}^{\star})}\right]$ in step 3, which leads to the final bound in step 4.

\noindent\textit{Step 1: Upper bound of $q\max_{j \in \mathcal{A}^{S}}u_{ij}\mathbf{x}_{ij}-V_{i}^{E}(\mathbf{x})$.} Based on the maximization objective, each demand would be recommended to exactly $\theta$ supplies when $\gamma=\theta$. Therefore, the objective function can be decomposed as a summation over these recommendations, from which an upper bound follows by bounding each term individually:
\begin{align*}
    q\max_{j \in \mathcal{A}^{S}}u_{ij}\mathbf{x}_{ij}-V_{i}^{E}(\mathbf{x}) & \le \left(\max_{j \in \mathcal{A}^{S}}u_{ij}\mathbf{x}_{ij}-a_i\right)\sum_{k=2}^{\theta}p(1-p)^{k-1} \\ 
   & \le \left(b_i-a_i\right)\sum_{k=2}^{\theta}p(1-p)^{k-1} \\  
    & = \left(b_i-a_i\right)\left(q-p\right).  
\end{align*}
\textit{Step 2: Upper bound of the objective difference $Obj^{E} (\mathbf{x}^{\star})-Obj^{E}(\mathbf{x}^{A})$.} Based on \hyperref[lemma:bounds]{Lemma \ref*{lemma:bounds}}, we can prove that
\begin{align*}
    Obj^{E} (\mathbf{x}^{\star})-Obj^{E}(\mathbf{x}^{A})  & = \sum_{i \in \mathcal{A}^{D}} V_{i}^{E}(\mathbf{x}^{\star})-\sum_{i \in \mathcal{A}^{D}} V_{i}^{E}(\mathbf{x}^{A}) \\
    & \le q\sum_{i \in \mathcal{A}^{D}}V_{i}^{A}(\mathbf{x}^{\star}) - \tau q|\mathcal{A}^{D}| \log(p)-\sum_{i \in \mathcal{A}^{D}} V_{i}^{E}(\mathbf{x}^{A})  \\
    & \le q\sum_{i \in \mathcal{A}^{D}}V_{i}^{A}(\mathbf{x}^{A}) - \tau q|\mathcal{A}^{D}| \log(p)-\sum_{i \in \mathcal{A}^{D}} V_{i}^{E}(\mathbf{x}^{A})  \\
    & \le \tau q|\mathcal{A}^{D}| \log(p)+ q\sum_{i \in \mathcal{A}^{D}}\max_{j \in \mathcal{A}^{S}}u_{ij}\mathbf{x}_{ij}^{A} +  \tau q|\mathcal{A}^{D}| \log (\theta) - \tau q|\mathcal{A}^{D}| \log(p)-\sum_{i \in \mathcal{A}^{D}} V_{i}^{E}(\mathbf{x}^{A})   \\
    & \le \tau q|\mathcal{A}^{D}| \log (\theta)+ \left(q-p \right)\sum_{i \in \mathcal{A}^{D}}\left(b_i-a_i\right),
\end{align*}
where the first inequality follows from the second and third inequalities of \hyperref[lemma:bounds]{Lemma \ref*{lemma:bounds}}, the second inequality is due to the optimality of $\mathbf{x}^{A}$, the third inequality follows from the fourth inequality of \hyperref[lemma:bounds]{Lemma \ref*{lemma:bounds}}, and the last inequality is from the result of step 1. Note that the above bound holds for any realization of $u_{ij}$, so we have $\mathbb{E}\left[Obj^{E} (\mathbf{x}^{\star})-Obj^{E}(\mathbf{x}^{A})\right] \leq \tau q|\mathcal{A}^{D}| \log (\theta)+ \left(q-p \right)\sum_{i \in \mathcal{A}^{D}}\left(b_i-a_i\right)$.

\noindent \textit{Step 3: Lower bound of  $\mathbb{E}\left[{Obj^{E} (\mathbf{x}^{\star})}\right]$.} We consider a solution $\mathbf{x}^{\square}$ that recommends demands to supplies following a prespecified order (e.g., recommending demand $i$ to supplies $\theta(i-1)+1,\dots, \theta\cdot i$). This solution is feasible for any realization of $u_{ij}$, and because the utility values are i.i.d., the utilities from the recommended supplies to demand $i$ still follow a uniform distribution over $[a_i,b_i]$. 
Let $Y_k^i$, $k=1,\dots,\theta$, denote the utility from the $k$th supplier that is recommended demand $i$. The number of supplies that accept the recommendation is denoted as $M=\sum_{k=1}^{\theta} \tilde{\xi}^i_{k}$ that follows a binomial distribution $\mathrm{Binomial}(\theta,p)$, where $\tilde{\xi}^i_k$ denotes the uncertain acceptance response from the $k$th supply that is recommended demand $i$. Conditioning on $M$, we obtain
\[\mathbb{E}\left[\max_{k} Y_k^i\tilde{\xi}^i_k\right] = \sum_{m=0}^{\theta} \Pr(M=m)\, \mathbb{E}\left[\max_{k} Y_k^i\tilde{\xi}^i_k\mid M=m\right].\]

If $m=0$, the maximum is zero. For $m\ge1$, independence and exchangeability imply that \[\mathbb{E}\left[\max_{k} Y_k^i\tilde{\xi}^i_k\mid M=m\right] = a_i + (b_i-a_i)\frac{m}{m+1}.\]

Hence, we have \begin{align*}
    \mathbb{E}\left[\max_{k} Y_k^i\tilde{\xi}^i_k\right] &= \sum_{m=1}^{\theta} \binom{\theta}{m} p^{m}(1-p)^{\theta-m}\Bigl[a_i+(b_i-a_i)\tfrac{m}{m+1}\Bigr]\\
    &=a_i\sum_{m=1}^{\theta}\binom{\theta}{m}p^{m}(1-p)^{\theta-m} + (b_i-a_i)\sum_{m=1}^{\theta}\binom{\theta}{m}p^{m}(1-p)^{\theta-m}\frac{m}{m+1}\\
    &=a_i(1-(1-p)^{\theta})+(b_i-a_i)\left(1-\frac{1-(1-p)^{\theta+1}}{(\theta+1)p}\right)
\end{align*}
Thus, $\mathbb{E}\left[Obj^{E} (\mathbf{x}^{\star})\right]$ satisfies
\begin{align*}
    \mathbb{E}\left[Obj^{E} (\mathbf{x}^{\star})\right] \ge \mathbb{E}\left[Obj^{E} (\mathbf{x}^{\square})\right] &= \sum_{i \in \mathcal{A}^{D}}\left(a_i(1-(1-p)^{\theta})+(b_i-a_i)\left(1-\frac{1-(1-p)^{\theta+1}}{(\theta+1)p}\right)\right)\\
    &=q\sum_{i \in \mathcal{A}^{D}}a_i+\left(\frac{\theta-\frac{q}{p}+q}{\theta+1}\right)\sum_{i \in \mathcal{A}^{D}}\left(b_i-a_i\right)
\end{align*}
\textit{Step 4: Bound of the approximation gap.} Combining the lower bound of the numerator (step 2) and the upper bound of the denominator (step 3), we derive the upper bound of the expected approximation gap as claimed. \Halmos
\subsection{Proof of Corollary \ref{cor1}}
\textit{Proof.} When the acceptance probabilities are homogeneous, the lower bound $\underline{p}$ and upper bound $\overline{p}$ coincide. In this case, based on the result of \hyperref[theoemhomo]{Theorem \ref*{theoemhomo}}, when $\theta = 1$, since $q=p$, each term $\tau |\mathcal{A}^{D}| \log (\theta)$, and $\left(1-\frac{p}{q} \right)\sum_{i \in \mathcal{A}^{D}}\left(b_i-a_i\right)$ in the numerator of the approximation performance bound becomes zero, thereby proving the corollary. \Halmos
\subsection{Proof of Proposition \ref{dominate}}
{\color{black}\textit{Proof.} Based on \hyperref[dapnew]{Proposition \ref*{dapnew}} and \hyperref[theoemhomo]{Theorem \ref*{theoemhomo}}, we can prove that
\begin{align*}
        & \mathbb{E}\left[Obj^{E} (\mathbf{x}^{\star})\right]-Obj^{E} (\mathbf{x}^{M}) \ge q\sum_{i \in \mathcal{A}^{D}}a_i+\left(\frac{\theta-\frac{q}{p}+q}{\theta+1}\right)\sum_{i \in \mathcal{A}^{D}}\left(b_i-a_i\right)-\sum_{i\in \mathcal{A}^D} p\cdot b_i,\\
     &  Obj^{E} (\mathbf{x}^{\star})-Obj^{E}(\mathbf{x}^{A}) \le \tau q|\mathcal{A}^{D}| \log (\theta)+ \left(q-p \right)\sum_{i \in \mathcal{A}^{D}}\left(b_i-a_i\right).
     \end{align*}
     Then, we have,
     \begin{align*}
     &\mathbb{E}\left[Obj^{E}(\mathbf{x}^{A})-Obj^{E} (\mathbf{x}^{M})\right] =\mathbb{E}\left[Obj^{E} (\mathbf{x}^{\star})\right]-Obj^{E} (\mathbf{x}^{M}) - \left(Obj^{E} (\mathbf{x}^{\star})-Obj^{E}(\mathbf{x}^{A})\right)\\
     &\ge q\sum_{i \in \mathcal{A}^{D}}a_i+\left(\frac{\theta-\frac{q}{p}+q}{\theta+1}\right)\sum_{i \in \mathcal{A}^{D}}\left(b_i-a_i\right)-p\sum_{i\in \mathcal{A}^D}b_i-\tau q|\mathcal{A}^{D}| \log (\theta)- \left(q-p \right)\sum_{i \in \mathcal{A}^{D}}\left(b_i-a_i\right)\\
     &=\sum_{i \in \mathcal{A}^{D}}a_i\left(q+\alpha\left(\frac{\theta-\frac{q}{p}+q}{\theta+1}\right)-p\left(1+\alpha\right)-\alpha\left(q-p\right)\right)-\tau q|\mathcal{A}^{D}| \log (\theta)\\
     & \approx \sum_{i \in \mathcal{A}^{D}}a_i\left(q-p+\alpha\left(\frac{\theta \left(1-q\right)-\frac{q}{p}}{\theta+1}\right)\right) \qquad \text{(when }\tau\to 0\text{)}
\end{align*}
Let $K(\theta,p)=\bigl(\theta(1-q)-q/p\bigr)/(\theta+1)$. To show $K(\theta,p)<0$, let $r=1-p\in(0,1)$. Then $q=1-r^\theta$ and
\[
\frac{q}{p}=\frac{1-r^\theta}{1-r}=\sum_{k=0}^{\theta-1} r^k.
\]
Hence
\[
\theta(1-q)-\frac{q}{p}=\theta r^\theta-\sum_{k=0}^{\theta-1}r^k\le \theta r^\theta-\theta r^{\theta-1}<0,
\]
where the inequality uses $\sum_{k=0}^{\theta-1}r^k\ge \theta r^{\theta-1}$ and the strict inequality uses $r^\theta<r^{\theta-1}$. Therefore $K(\theta,p)<0$.
Since $\sum_{i \in \mathcal{A}^{D}}a_i >0$, $\mathbb{E}\left[Obj^{E}(\mathbf{x}^{A})-Obj^{E} (\mathbf{x}^{M})\right]>0$ is equivalent to $(q-p)+\alpha K(\theta,p)>0$. i) If $\theta=1$, then $q=p$, so $\alpha K(\theta,p)<0$ for all $\alpha>0$. ii) If $\theta\ge 2$, then $q-p>0$ and $K(\theta,p)<0$, so the inequality holds if and only if
\[
\alpha<-\frac{q-p}{K(\theta,p)}
=\frac{(q-p)(\theta+1)}{\,q/p-\theta(1-p)^\theta\,}
=\alpha_{\max}(\theta,p).
\]
Moreover, $\alpha_{\max}(\theta,p)>0$ since $q-p>0$ and $q/p-\theta(1-p)^\theta>0$ for all $p\in(0,1)$, thereby proving the proposition.\halmos}
\subsection{Proof of Theorem \ref{theoremheter}}
\textit{Proof.} Similar to \hyperref[theoemhomo]{Theorem \ref*{theoemhomo}}, we prove the result in the following steps. 

\noindent\textit{Step 1: Lower bound of $V_{i}^{A}(x)$ in terms of $V_{i}^{E}(x)$.} Substituting the upper bound of $V_{i}^{E}(\mathbf{x})$ into the lower bound of $V_{i}^{A}(\mathbf{x})$ by \hyperref[lemma:bounds]{Lemma \ref*{lemma:bounds}}, we have:
$$ V_{i}^{A}(\mathbf{x}) \ge \tau \log(\underline{p})+ U_{i}(\mathbf{x}) \ge \frac{V_{i}^{E}(\mathbf{x})}{\overline{q}} + \tau \log(\underline{p}). \quad \text{(P1)}$$

\noindent \textit{Step 2: Upper bound of $V_{i}^{A}(x)$ in terms of $V_{i}^{E}(x)$.} Substituting the lower bound of $V_{i}^{E}(\mathbf{x})$ into the upper bound of $V_{i}^{A}(\mathbf{x})$ by \hyperref[lemma:bounds]{Lemma \ref*{lemma:bounds}}, we have:
$$ V_{i}^{A}(\mathbf{x}) \le \frac{V_{i}^{E}(\mathbf{x})}{\underline{p}} + \tau \log(\theta)+\tau \log(\overline{p}). \quad \text{(P2)}$$

\noindent \textit{Step 3: Relationship between $Obj^{E} (\mathbf{x}^{A})$ and $Obj^{E}(x^{\star}$).} Since $x^{A}$ is an optimal solution for A-SP, we have $Obj^{A} (\mathbf{x}^{A})\geq Obj^{A} (\mathbf{x}^{\star})$. Summing the inequalities over all $i\in \mathcal{A}^{D}$:
 \begin{align*}
     &Obj^{A} (\mathbf{x}^{A})=\sum_{i\in \mathcal{A}^{D}}V_{i}^{A}(\mathbf{x}^{A})\le \frac{Obj^{E} (\mathbf{x}^{A})}{\underline{p}} + \tau  \left| \mathcal{A}^{D}\right|\log(\theta)+\tau  \left| \mathcal{A}^{D}\right| \log(\overline{p}) & \text{Applying (P2) to } \mathbf{x}^{A}\\
     &Obj^{A} (\mathbf{x}^{\star})=\sum_{i\in \mathcal{A}^{D}}V_{i}^{A}(\mathbf{x}^{\star})\ge \frac{Obj^{E}(\mathbf{x}^{\star})}{\overline{q}}+\tau  \left| \mathcal{A}^{D}\right| \log(\underline{p}) & \text{Applying (P1) to } \mathbf{x}^{\star}
 \end{align*}
 Combining these with the optimality of $x^{A}$ for A-SP, we have
 $$\frac{Obj^{E} (\mathbf{x}^{A})}{\underline{p}} + \tau  \left| \mathcal{A}^{D}\right|\log(\theta)+\tau  \left| \mathcal{A}^{D}\right| \log(\overline{p})\ge \frac{Obj^{E}(\mathbf{x}^{\star})}{\overline{q}}+\tau  \left| \mathcal{A}^{D}\right| \log(\underline{p})$$
 which is equivalent to
 $$Obj^{E} (\mathbf{x}^{A})\ge \frac{\underline{p}}{\overline{q}}Obj^{E}(\mathbf{x}^{\star})-\tau\underline{p}\left| \mathcal{A}^{D}\right| \log(\frac{\theta\overline{p}}{\underline{p}})$$
 
\noindent \textit{Step 4: Lower bound of $Obj^{E} (\mathbf{x}^{\star})$.} Consider an order $i$ that is recommended to $m$ suppliers and that the utility of any supplier is greater than $a_i$, we have $V_{i}^{E}(\mathbf{x}) \ge \mathbb{E}_{\mathbb{P}} \left [\max_{j\in \mathcal{A}^{S}} a_{i} \tilde{\xi}_{ij} x_{ij}\right] = \left( 1- \Pi_{j\in \mathcal{A}^{S},x_{ij}=1}(1-p_{ij}) \right)a_{i}  \ge\left(1-(1-\underline{p})^{m}\right)a_{i}$, where we leverage the independence of supplier responses facing the same demand. Moreover, i) when $\gamma<\theta$, there exists a feasible solution in which each demand is first recommended to $\lfloor \gamma \rfloor$ suppliers, and subsequently, if set $|\mathcal{A}^{D\prime}|$ is nonempty, the demands in $|\mathcal{A}^{D\prime}|$ are further recommended to one additional supplier, then, we have
\begin{align*}
    Obj^{E} (\mathbf{x}^{\star}) & \ge 
   \sum_{i\in \mathcal{A}^{D\prime}}\left(1-(1-\underline{p})^{\lfloor \gamma \rfloor+1}\right)a_{i}+ \sum_{i\in \mathcal{A}^{D}\setminus \mathcal{A}^{D\prime}}\left(1-(1-\underline{p})^{\lfloor \gamma \rfloor}\right)a_{i}
\end{align*}

\textcolor{black}{ii) When $\gamma\ge \theta$, similarly, we have 
$$ Obj^{E} (\mathbf{x}^{\star}) \ge \sum_{i\in \mathcal{A}^{D}}\left(1-(1-\underline{p})^{\theta}\right)a_{i}.$$
By revising the definition of the set $\mathcal{A}^{D\prime}$ so that it becomes empty when $\gamma\ge \theta$, we unify the two cases and obtain:
$$Obj^{E} (\mathbf{x}^{\star})  \ge 
   \sum_{i\in \mathcal{A}^{D\prime}}\left(1-(1-\underline{p})^{\lfloor \gamma \rfloor+1}\right)a_{i}+ \sum_{i\in \mathcal{A}^{D}\setminus \mathcal{A}^{D\prime}}\left(1-(1-\underline{p})^{\min(\lfloor \gamma \rfloor, \theta)}\right)a_{i}$$}
\noindent \textit{Step 5: Bound of the approximation gap.} Finally, combining the above results leads to
\begin{align*}&\frac{Obj^{E} (\mathbf{x}^{\star})-Obj^{E} (\mathbf{x}^{A})}{Obj^{E} (\mathbf{x}^{\star})}\le 1-\frac{\underline{p}}{\overline{q}}+\frac{\tau\underline{p}\left| \mathcal{A}^{D}\right| \log(\frac{\theta\overline{p}}{\underline{p}})}{Obj^{E} (\mathbf{x}^{\star})}\\
 &\le  1-\frac{\underline{p}}{\overline{q}}+\frac{\tau\underline{p}\left| \mathcal{A}^{D}\right| \log(\frac{\theta\overline{p}}{\underline{p}})}{\sum_{i\in \mathcal{A}^{D\prime}}\left(1-(1-\underline{p})^{\lfloor \gamma \rfloor+1}\right)a_{i}+ \sum_{i\in \mathcal{A}^{D}\setminus \mathcal{A}^{D\prime}}\left(1-(1-\underline{p})^{\min(\lfloor \gamma \rfloor, \theta)}\right)a_{i}}.\halmos\end{align*} 

\subsection{Proof of Corollary \ref{corogeneral}}
\textit{Proof.} When there exists correlation, the upper bound of $V_{i}^{E}(\mathbf{x})$ in \hyperref[lemma:bounds]{Lemma \ref*{lemma:bounds}} would take the form of $ V_i^{E}(\mathbf{x}) \le U_{i}(\mathbf{x})$ and the lower bound of $Obj^{E} (\mathbf{x}^{\star})$ in \hyperref[theoremheter]{Theorem \ref*{theoremheter}} would be $Obj^{E} (\mathbf{x}^{\star})  \ge \sum_{i\in \mathcal{A}^{\underline{D}}}\underline{p}a_{i}$. The proof is similar to \hyperref[theoremheter]{Theorem \ref*{theoremheter}}, and we have 
\begin{align*}
    \frac{Obj^{E} (\mathbf{x}^{\star})-Obj^{E} (\mathbf{x}^{A})}{Obj^{E} (\mathbf{x}^{\star})} &\le \frac{Obj^{E} (\mathbf{x}^{\star})-\left(\underline{p}Obj^{E}(\mathbf{x}^{\star})-\tau\underline{p}\left| \mathcal{A}^{D}\right| \log(\frac{\theta\overline{p}}{\underline{p}})\right)}{Obj^{E} (\mathbf{x}^{\star})} \\
    & \le 1-\underline{p}+\frac{\tau\left| \mathcal{A}^{D}\right| \log(\frac{\theta\overline{p}}{\underline{p}})}{\sum_{i\in \mathcal{A}^{\underline{D}}}a_{i}}.\halmos
\end{align*}

\subsection{Proof of Proposition \ref{propsub}}
{\color{black} \textit{Proof.} 
For a feasible recommendation set \((S_i)_{i\in\mathcal A^D}\), we denote the objective as $\sum_{i\in\mathcal A^D} F_i(S_i)$, we prove the proposition in four steps.

\noindent \textit{Step 1: Structural property of the objective.} By Proposition 5.2 in \citet{ekbatani2026lyft1}, for every \(i\in\mathcal A^D\), the objective function \(F_i(\cdot)\) is monotone and submodular. Hence the total objective
\[
Obj^E(x):=\sum_{i\in\mathcal A^D}F_i(S_i(x)),
\qquad
S_i(x):=\{j\in\mathcal A^S:x_{ij}=1\},
\]
is a monotone submodular set function in the recommendation pairs \((i,j)\). 

\noindent \textit{Step 2: Continuous-greedy algorithm with the additional cap \(\theta\).} We introduce the extra per-demand recommendation cap to the submodular welfare maximization problem of \citet{ekbatani2026lyft1}. The fractional feasible region is derived by relaxing the binary variables to \(x_{ij}\in[0,1]\) and considering the polytope
\[
P_\theta:=\left\{x\in[0,1]^{|\mathcal A^D|\times |\mathcal A^S|}:
\sum_{i\in\mathcal A^D}x_{ij}\le 1,\ \forall j\in\mathcal A^S,
\quad
\sum_{j\in\mathcal A^S}x_{ij}\le \theta,\ \forall i\in\mathcal A^D
\right\}.
\]
This polytope is down-monotone: if \(0\le y\le x\) coordinatewise and \(x\in P_\theta\), then \(y\in P_\theta\) as well. Moreover, linear optimization over
\(P_\theta\) is a bipartite \(b\)-matching LP with capacity \(1\) on each supplier and capacity \(\theta\) on each demand, and is therefore solvable in polynomial time. Let \(\bar f(x)\) be the multilinear extension of the discrete objective. Since \(\bar f\) is smooth monotone submodular, and \(P_\theta\) is a down-monotone polytope, the standard continuous-greedy guarantee applies \citep{vondrak2008optimal}, leading to a $(1-1/e)$-approximation in the idealized continuous-time process. For any prescribed $\varepsilon>0$, by discretizing the continuous trajectory sufficiently finely and estimating marginal values to sufficient accuracy,  the resulting finite algorithm incurs at most an \(\varepsilon\) loss. Consequently, it yields a fractional point \(x^{\mathrm{CG}}\in P_\theta\) such that
\[
\mathbb{E}[\bar f(x^{\mathrm{CG}})]
\ge (1-1/e-\varepsilon)\bar f(\mathbf{x}^\star)=
(1-1/e-\varepsilon)Obj^{E} (\mathbf{x}^{\star}).
\]

\noindent \textit{Step 3: Rounding with a finite number of steps yields an integer feasible solution.}
To transform the fractional solution \(x^{\mathrm{CG}}\) into an integral feasible solution without decreasing the multilinear value, we construct the fractional support graph
$$G_x=(\mathcal A^D\cup\mathcal A^S,E_x),
\qquad
E_x:=\{(i,j):0<x_{ij}<1\}.$$
As long as \(E_x\neq\varnothing\), select one connected component of \(G_x\) containing at least one edge. This component necessarily contains either a simple cycle or, if acyclic, a maximal simple path.
If the component contains a cycle, choose a simple cycle and assign alternating signs \(+1,-1,+1,-1,\dots\) to its consecutive edges. If the component is acyclic, choose a maximal simple path and again assign alternating signs \(+1,-1,+1,-1,\dots\) along the path. Fix a fractional point \(x\in P_\theta\), and let \(d\) be the alternating direction induced by the selected maximal simple path or cycle in the support graph of \(x\). Define
\[
\mathcal I(x,d):=\{\lambda\in\mathbb R: x+\lambda d\in P_\theta\}.
\]
Since \(x\in P_\theta\), we have \(0\in\mathcal I(x,d)\), so \(\mathcal I(x,d)\neq\varnothing\). Moreover, because \(P_\theta\) is described by linear capacity constraints together with box constraints, \(\mathcal I(x,d)\) is a closed interval, which we write as
\[
\mathcal I(x,d)=[-\alpha,\beta]
\]
for some \(\alpha,\beta\ge 0\) to ensuring that each variable lies in $[0,1]$.

By construction of the alternating direction, along a selected cycle, the net change at every visited demand node and supplier node is zero; along a selected path, the same holds for every internal node. Furthermore, because the selected path is maximal in the fractional graph, its two endpoints have exactly one incident fractional edge in $G_x$. Since the capacities $1$ and $\theta$ are integral, the capacity constraints at these endpoints cannot be tight (i.e., an integer capacity cannot be perfectly met by a sum of integers and exactly one strictly fractional value). Hence, feasibility with respect to all capacity constraints is preserved along the segment \(x+\lambda d\), and the endpoints
\[
x^-:=x-\alpha d,\qquad x^+:=x+\beta d
\]
are both feasible. By maximality of \(\alpha\) and \(\beta\), at least one previously fractional coordinate becomes \(0\) or \(1\) at each endpoint. Therefore, replacing \(x\) by either \(x^-\) or \(x^+\) strictly decreases the number of fractional coordinates. Since the number of fractional coordinates is finite and each iteration strictly decreases this number, the procedure terminates after finite steps at an integral point in \(P_\theta\).

\noindent \textit{Step 4: Lossless rounding explanation.}
For the alternating direction \(d\) constructed in Step 3, we define
$\phi(\lambda):=\bar f(x+\lambda d),\ \lambda\in[-\alpha,\beta]$.
We claim that \(\phi(\lambda)\) is convex on \([-\alpha,\beta]\). Indeed, since the objective is separable across demands, we can write $\bar f(x)=\sum_{i\in\mathcal A^D}\bar f_i(x_i)$, where 
$x_i:=(x_{ij})_{j\in\mathcal A^S}$. Hence, we have 
$\phi(\lambda)=\sum_{i\in\mathcal A^D}\bar f_i(x_i+\lambda d_i)$,
where \(d_i:=(d_{ij})_{j\in\mathcal A^S}\). By construction of the alternating path/cycle, for each fixed demand \(i\), the vector \(d_i\) has at most two nonzero entries. Therefore, \(\bar f_i(x_i+\lambda d_i)\) is either constant, affine, or of the form $\bar f_i(\ldots,x_{ij}+\lambda,\ldots,x_{ik}-\lambda,\ldots)$
for some \(j,k\in\mathcal A^S\). In the first two cases, the function is trivially convex. In the third case, differentiating twice yields
\[
\frac{d^2}{d\lambda^2}\bar f_i(\ldots,x_{ij}+\lambda,\ldots,x_{ik}-\lambda,\ldots)=\frac{\partial^2\bar f_i}{\partial x_{ij}^2}+
\frac{\partial^2\bar f_i}{\partial x_{ik}^2}-
2\frac{\partial^2\bar f_i}{\partial x_{ij}\partial x_{ik}}.
\]

Since \(\bar f_i\) is the multilinear extension of \(F_i\), it is affine in each individual coordinate when all other coordinates are fixed, hence $\frac{\partial^2 \bar f_i}{\partial x_{ij}^2}=\frac{\partial^2 \bar f_i}{\partial x_{ik}^2}=0$. Moreover, since \(F_i\) is submodular, the multilinear extension satisfies $\frac{\partial^2 \bar f_i}{\partial x_{ij}\partial x_{ik}}\le 0$. Therefore, the second derivative of $\bar f_i(x_i+\lambda d_i)$ is nonnegative, and \(\bar f_i(x_i+\lambda d_i)\) is convex in \(\lambda\). Summing over all \(i\in\mathcal A^D\), we conclude that \(\phi(\lambda)\) is convex on \([-\alpha,\beta]\). Because \(0\in[-\alpha,\beta]\), convexity implies
$\phi(0)\le \max\{\phi(-\alpha),\phi(\beta)\}$,
that is, $\bar f(x)\le \max\{\bar f(x-\alpha d),\,\bar f(x+\beta d)\}$. Therefore, replacing \(x\) by the better endpoint never decreases the multilinear objective. Combining this fact with Step 3, the procedure terminates after finite iterations at an integral feasible point \(\widehat x\) satisfying
\[
\bar f(\widehat x)\ge \bar f(x^{\mathrm{CG}}).
\]
Finally, since \(\widehat x\) is integral, the multilinear extension coincides with the original objective. Therefore,
\[
\mathbb{E}[Obj^{E} (\mathbf{x}^{S})]=\mathbb{E}[\bar f(\widehat x)]\ge\mathbb{E}[\bar f(x^{\mathrm{CG}})]\ge
(1-1/e-\varepsilon)Obj^{E} (\mathbf{x}^{\star}).
\]
This completes the proof.
\Halmos}

\section{Additional Algorithm and Model Details}
\subsection{Cap-SW algorithm pseudocode} \label{sec:sub}
\begin{algorithm}[t]
\caption{\textcolor{black}{Cap-SW: Monte Carlo Continuous Greedy with Rounding}}
\label{alg:ba-cg-theta}
\begin{algorithmic}[1]
\Require $\mathcal{A}^{D}$, $\mathcal{A}^{S}$, $\{u_{ij}\}$, $\mathbb{P}$, $\theta$
\Ensure An integral feasible recommendation decision $x=\{x_{ij}\}$

\State Compute the maximum number of iterations $T$
\State Set $\Delta \gets 1/T$, define
$f_i(S):=\mathbb{E}_{\mathbb{P}}\!\left[\max_{j\in S} u_{ij}\tilde{\xi}_{ij}\right]$,
initialize $x_{ij}\gets 0$ for all $i\in\mathcal{A}^{D},j\in\mathcal{A}^{S}$

\For{$t=0,1,\ldots,T-1$}

    \Statex \hspace{\algorithmicindent}// \textit{Monte Carlo gradient estimation}
    \ForAll{$i\in\mathcal{A}^{D},\,j\in\mathcal{A}^{S}$}
        \State Generate Monte Carlo samples of a random supplier subset $R_{x,i}\subseteq \mathcal A^S$, where each supplier $j\in\mathcal A^S$ is included independently with probability $x_{ij}$
        \State Estimate the marginal gain of edge $(i,j)$ by $\hat g_{ij}\gets \widehat{\mathbb E}\left[
        f_i(R_{x,i}\cup\{j\})-f_i(R_{x,i}\setminus\{j\})
        \right]$
        
    \EndFor

    \Statex \hspace{\algorithmicindent}// \textit{Linear oracle}
    \State Compute a search direction $v=\{v_{ij}\}$ by solving
$v\in\arg\max_{v\in P_\theta}
    \sum_{i\in\mathcal{A}^{D}}\sum_{j\in\mathcal{A}^{S}}\hat g_{ij}v_{ij}$,
    where
    \[
    P_\theta:=\left\{x\in[0,1]^{|\mathcal{A}^{D}|\times|\mathcal{A}^{S}|}:
    \sum_{j\in\mathcal{A}^{S}}x_{ij}\le \theta,\ \forall i\in\mathcal{A}^{D},\ 
    \sum_{i\in\mathcal{A}^{D}}x_{ij}\le 1,\ \forall j\in\mathcal{A}^{S}\right\}.
    \]

    \State Update $x\gets x+\Delta v$
    \State Estimate the current objective value by Monte Carlo sampling 
\EndFor

\Statex // \textit{Rounding}
\While{$x$ is fractional}
    \State Construct the fractional support graph
    $G_x=(\mathcal{A}^{D}\cup\mathcal{A}^{S},E_x),\
    E_x:=\{(i,j):0<x_{ij}<1\}$
    \State Find an alternating single path or cycle $C$ in $G_x$ and define the alternating direction $d$ on the edges of $C$
    \State Compute the maximal feasible interval
    $[-\alpha,\beta]:=\{\lambda\in\mathbb R:x+\lambda d\in P_\theta\}$
    \State Let $x^{-}\gets x-\alpha d$ and $x^{+}\gets x+\beta d$
    \State Set $x\gets x^{+}$ if $\bar f(x^{+})\ge \bar f(x^{-})$; otherwise set $x\gets x^{-}$
\EndWhile
\State \Return $x$
\end{algorithmic}
\end{algorithm}

\subsection{Performance evaluation of L-shaped method}\label{addLshaped}
\textcolor{black}{We implement the L-shaped method \citep{van1969shaped} for SP-SAA model, in which the master problem determines the recommendation decision $\mathbf{x}$, while the subproblems correspond to a series of matching problems under different scenarios, associated with the decision $\mathbf{\omega}_s$. For any feasible master-problem solution $\mathbf{\bar{x}}$, each subproblem is feasible and bounded. Therefore, throughout the L-shaped procedure, it suffices to add only the optimality cuts derived from duality. Specifically, let $\mathbf{\alpha^{s}}=(\alpha_{i}^{s})_{i \in \mathcal{A}^{D}}$ and $\mathbf{\beta^{s}}=(\beta_{ij}^{s})_{i \in \mathcal{A}^{D},j \in \mathcal{A}^{S}}$ denote the dual variables associated with constraints \hyperref[eq:b2]{(\ref*{eq:b2})} and \hyperref[eq:b3]{(\ref*{eq:b3})}, respectively. Then, the optimality cut for scenario $s$ can be written as $\eta_s \le \sum_{i \in \mathcal{A}^{D}}\bar{\alpha}^{s}+\sum_{i \in \mathcal{A}^{D}}\sum_{j\in \mathcal{A}^{S}}\xi_{ij}^{s}\bar{\beta}_{ij}^{s}x_{ij}$. To further strengthen the L-shaped method, we introduce the following acceleration strategies. i) Since the number of samples $|\mathcal{S}|$ is large, we aggregate the cuts from all scenarios into a single cut. ii) For each scenario, we construct a closed-form dual extreme-point cut. Let $u_{i,max}^s$ be the maximum utility among all active edges, i.e., $\max\{u_{ij}:\xi_{ij}^s\bar{x}_{ij}>0\}$. If it does not exist, then set $u_{i,max}^s=0$. Then, the corresponding dual solution is given by $\alpha_{i}^{s}=u_{i,max}^s$ and $\beta_{ij}^{s}=\max\{0,u_{ij}-\alpha_i^s\}$. iii) We implement the algorithm within a branch-and-cut framework since the master problem is relatively difficult to solve. iv) We incorporate warm-start heuristics and branching priorities. However, the computational performance of the L-shaped method is inferior to that of directly solving the SAA model by Gurobi.}

\textcolor{black}{Tables \ref{tab:moderate} and \ref{tab:large} report the results for the moderate-scale and large-scale instances, including the average optimality gap reported by the solver (in columns \textit{GAP (\%)}) and the number of instances for which the out-of-sample expected utility of L-shaped method is higher than that of SAA (in columns \textit{L-shaped-Impr}). Overall, we do not observe significant computational improvements from the implemented L-shaped decomposition algorithm: although the L-shaped method can improve on some instances, the SAA implementation yields generally better solutions on both moderate and large instances, and the gap between the two algorithms narrows as problem size increases. We believe this may be due to Gurobi's enhanced solution capabilities for large-scale MILPs. Similar results have been documented in the related literature: directly solving the SAA by CPLEX can be more efficient than a tailored decomposition method, as shown in Table 3 of \citeapp{oliveira2014accelerating} and Table 4 of \citeapp{rebennack2020two}.}

\begin{table}[htbp]
        \caption{Performance evaluation for moderate-scale instances}
        \label{tab:moderate}
        \centering
        \begin{adjustbox}{width=0.6\textwidth}
\begin{tabular}{cccc}
\hline
Sample size & SAA-GAP(\%) & L-shaped-GAP(\%) & L-shaped-Impr \\ \hline
100         & 0.53        & 7.22             & 0/10          \\
500         & 5.71        & 9.58             & 3/10          \\
1000        & 7.97        & 9.70             & 4/10          \\
2000        & 9.30        & 11.81            & 3/10          \\ \hline
\end{tabular}
\end{adjustbox}
    \end{table}
    \begin{table}[htbp]
    
        \caption{Performance evaluation for large-scale instances}
        \label{tab:large}
        \centering
        \begin{adjustbox}{width=0.6\textwidth}
\begin{tabular}{cccc}
\hline
Sample size & SAA-GAP(\%) & L-shaped-GAP(\%) & L-shaped-Impr \\ \hline
100         & 9.266       & 23.04            & 0/10          \\
150         & 21.943      & 23.73            & 4/10          \\ \hline
\end{tabular}
\end{adjustbox}
    \end{table}
\subsection{Formulation of NPP}\label{sec:anpp}
The objective of the NPP model is to minimize the total distance between drivers’ locations and order origins $\sum_{i \in \mathcal{A}^{D}} \sum_{j \in \mathcal{A}^{S}} d_{ij}x_{ij}$ (while satisfying the vehicle type constraints). In our utility maximization formulation, we take the utility $u_{ij} = M^{D} - d_{ij}$ for a large $M^D$ ($M^D$ is needed to avoid trivial solutions of setting $x_{ij}$'s to be zero):
    \begin{align}
       \textbf{[NPP]}:  \max \ & \sum_{i \in \mathcal{A}^{D}} \sum_{j \in \mathcal{A}^{S}} (M^{D}-d_{ij})x_{ij} \\
       {\rm s.t.} \ & \hyperref[eq:a1]{(\ref*{eq:a1})-(\ref*{eq:a3})} \nonumber.
    \end{align}
\bibliographystyleapp{informs2014trsc}
\bibliographyapp{appealref.bib}

\end{document}